\documentclass{amsart}

\usepackage{mathrsfs}
\usepackage{epic}
\usepackage{pstricks}

\newtheorem{theorem}{Theorem}[section]

\theoremstyle{definition}
\newtheorem{definition}[theorem]{Definition}

\newtheorem{corollary}[theorem]{Corollary}

\theoremstyle{remark}
\newtheorem{remark}[theorem]{Remark}

\numberwithin{equation}{section}

\begin{document}

\title{A survey on $q$-polynomials and
their orthogonality properties}

%    Information for first author
\author{Roberto S. Costas-Santos}
%    Address of record for the research reported here
\address{Department of Mathematics,
University of California, Santa
Barbara, California 93106, US}
%    Current address
%\curraddr{Department of Mathematics and Statistics,
%Case Western Reserve University, Cleveland, Ohio 43403}
\email{rscosa@gmail.com}
%    \thanks will become a 1st page footnote.
\thanks{RSCS acknowledges financial
support from the Ministerio de Ciencia e
Innovaci\'on of Spain, grant
MTM2009-12740-C03-01, and the
program of postdoctoral grants
(Programa de becas postdoctorales).}

%    Information for second author
\author{Joaquin F. S\'{a}nchez-Lara}l
\address{Departamento de Matem\'{a}tica
Aplicada, Facultad de CC. Econ\'omicas
y Empresariales, Universidad de Granada,
Campus de la cartuja, s/n. 18071 Granada,
Spain}
\email{jslara@ugr.es}
\thanks{JFSL acknowledges financial
support from the Spanish Ministry of
Education e Innovaci\'on, grant MTM2008-06689-C02,
and Junta de Andaluc\'{\i}a, grant FQM229.}

%    General info
\date{\today}

\begin{abstract}
In this paper we study the orthogonality
conditions satisfied by the classical
$q$-orthogonal polynomials
that are located at the top of the $q$-Hahn
tableau (big $q$-jacobi polynomials (b$q$J))
and the Nikiforov--Uvarov tableau (Askey-Wilson
polynomials (AW)) for almost any complex value
of the parameters and for all non-negative
integers degrees.

We state the degenerate version of Favard's
theorem, which is one of the keys of the paper,
that allow us to extend the orthogonality
properties valid up to some integer degree $N$
to Sobolev type orthogonality properties.

We also present, following an analogous process
that applied in \cite{cola1}, tables with the
factorization and the discrete Sobolev-type orthogonality
property for those families which satisfy a
finite orthogonality property, i.e. it consists
in sum of finite number of masspoints, such as
$q$-Racah ($q$R), $q$-Hahn ($q$H), dual $q$-Hahn (d$q$H),
and $q$-Krawtchouk polynomials ($q$K), among others.
\end{abstract}
\maketitle
%%%%%%%%%%%%%%%%%%%%%%%%%%%%%%%%%%%%%%%%%%%%%%%%%
\section{Introduction} \label{sec1}
The classical orthogonal polynomials
constitute a very important and
interesting set of special functions
and more specifically of orthogonal
polynomials. They are very
interesting mathematical objects
that have attracted the attention not
only of mathematicians since their
appearance at the end of the XVIII
century connected with some physical
problems.
They are used in several branches of
mathematical and physical sciences and
they have a lot of useful properties:
they satisfy a three-term recurrence
relation (TTRR), they are the solution
of a second order linear differential
(or difference) equation, their
derivatives (or finite differences)
also constitute an orthogonal family,
among others (for a recent review see
e.g. \cite{alv3}).

In this survey we are going to focus
on classical $q$-orthogonal polynomials
--also called $q$-polynomials-- which
are polynomial eigenfunctions of the
second order hypergeometric-type
homogeneous linear difference operator
\begin{equation} \label{1:1}
{\mathscr H}=\sigma(s)\frac{\Delta}
{\Delta x(s-\frac 12)}\frac{\nabla}
{\nabla x(s)}+\tau(x(s))\frac{\Delta}
{\Delta x(s)},
\end{equation}
where $\widehat \sigma(x(s))\stackrel{\rm
def}=\sigma(s)+\frac 12\tau(x(s))\Delta
x(s-\frac 12)$ and $\tau(x(s))$ are polynomials
on $x(s)$ with $\deg\widehat\sigma\le
2$ and $\deg \tau=1$.

In fact they appear in several branches
of the natural sciences, e.g., quantum
groups and algebras \cite{koo1,koo2,vikl},
quantum optics ,  continued fractions, theta
functions, elliptic functions, \dots; among
others \cite{and,fin,gara,nisuuv}.

$q$-polynomials have been intensively studied
by the American School starting from the
works of G. E. Andrews and R. Askey \cite{anas}
arising the $q$-Askey tableau, and the Russian
(former Soviet) school, starting from the works
in \cite{niuv3} and further developed by N. M.
Atakishiyev and S. K. Suslov (see
\cite{atrasu,atsu5,nisuuv,niuv2} and references
therein) arising the Nikiforov--Uvarov tableau.

It is known that any family of polynomial
eigenfunctions $(p_n)$ of \eqref{1:1} satisfies
a TTRR \cite{niuv3}, i.e. there exist two
sequences of complex numbers, $(\beta_n)$
and $(\gamma_n)$, such that for $n\ge 1$
\begin{equation} \label{1:2}
p_{n+1}(x)=(x-\beta_{n})p_{n}(x)
-\gamma_{n}p_{n-1}(x),
\end{equation}
with initial conditions $p_0(x)=1$,
$p_{1}(x)=x-\beta_0$.

On the other side, if a sequence of monic
polynomials $(p_n)\ $satisfies the initial
conditions $p_0(x)=1$, $p_{1}(x)=x-\beta_0$,
and the TTRR  \eqref{1:2} then these
polynomials are orthogonal with respect to the
moment functional ${\mathscr L}_0$
\cite[\S 1]{chi1}, defined by
${\mathscr L}_0(p_n)=\delta_{n,0}$,
$n\ge 0$, i.e., for $n\ne m$
\begin{equation}\label{1.3}
{\mathscr L}_0(p_np_m)=0.
\end{equation}
If $\gamma_n\neq 0$ for $n\geq 1$ then the
polynomials defined by \eqref{1:2}, are the
unique normal and monic polynomials satisfying
the orthogonality property \eqref{1.3}.
Moreover, if $\beta_n$, $\gamma_n$ are real,
$\gamma_n>0$, then there exists a
positive Borel measure $\mu$ such that
$$
{\mathscr L}_0(p)=\int_{\mathbb R}p d\mu.
$$
This result is known as Favard's theorem (see
\cite{fav}, \cite[p. 21]{chi1}), although
this result was also discovered (independent
of J. Favard) by I. P. Natanson in 1935
\cite{nat} and was presented by himself in a
seminar led by S. N. Bernstein. He then did not
publish the result since the work of
J. Favard appeared in the meantime.

Our main aims here are to study the orthogonality
conditions satisfied by Askey-Wilson and big
$q$-Jacobi polynomials for almost any complex
value of the parameters, any complex value
of $q$ and all non-negative integer degrees.
In all the cases, the proposed orthogonality
conditions characterizes such polynomials.
When there exists a $\gamma_N=0$ in \eqref{1:2},
an extension of the Favard's result is used in
order to establish a Sobolev-type orthogonality.
In such a case we also give the factorization
$$
p_{n+N}=p_N p_n^{(N)}, \qquad n\geq 0,
$$
where $p_n^{(N)}$ is the associated polynomial
of order $N$ and degree $n$, which also belongs
to the Nikiforov--Ubarov and/or $q$-Askey tableaux.
We present a table with the Sobolev-type
orthogonality and the factorization for all
the $q$-polynomials considered in Section
\ref{sec2.2} for whose TTRR there exists,
at least, one $N$ such that $\gamma_N=0$.

The structure of the paper is as follows.
The preliminaries which will be used
throughout the paper as well as the
extension of the Favard's theorem are
given in Section 2.
In Sections 3 and  4 we study the
orthogonality conditions for Askey-Wilson
and big $q$-Jacobi polynomials respectively
for almost any value of the complex
parameters.
In section 5 we study the orthogonality
conditions for Askey-Wilson for $|q|=1$.
In Section 6 we give a table for all the
families of $q$-polynomials which satisfy
a discrete orthogonality with a finite
number of masses like $q$-Racah, $q$-Hahn,
dual $q$-Hahn, and  $q$-Krawtchouk
polynomials, among others; and we finish
this paper giving some conclusions and
outlooks. An appendix is also included.
%%%%%%%%%%%%%%%%%%%%%%%%%%%%%%%%%%%%%%%%%%%%%%%%%
\section{Preliminaries} \label{sec2}
In this subsection we summarize some
definitions and preliminary results
that will be useful throughout the work.
Most of them can be found in \cite{chi1}.

\begin{definition} \label{def2.1}
Let $(\mu_n)$ be a sequence of complex numbers
(moment sequence) and $\mathscr L$ a
functional acting on the linear space
of polynomials $\mathbb P$ with complex
coefficients.
We say that $\mathscr L$ is a {\sf moment
functional associated with $(\mu_n)$} if
$\mathscr L$ is linear and
${\mathscr L}(x^n)=\mu_n$, $n\ge 0$.
\end{definition}
\begin{definition} \label{def2.2}
The polynomial sequence $(p_n)$ is a {\sf
orthogonal polynomials system} (OPS) with
respect to a moment functional $\mathscr L$
if the following conditions hold:
\begin{enumerate}
\item $p_n$ is a polynomial of exact
degree $n$, i.e. the polynomial sequence
$(p_n)$ is {\sf normal}.
\item ${\mathscr L}(p_np_m)=0$, $m\ne n$.
\item ${\mathscr L}(p_n^2)\ne 0$.
\end{enumerate}
\end{definition}
This third condition is imposed in order
to have a unique OPS:
if $\mathscr{L}(p_N^2)=0$ then
$$
\mathscr{L}((p_{N+1}+\alpha p_N) x^m)=0,\qquad
m=0,\dots,N,\qquad \forall \alpha\in \mathbb C.
$$
The next result is a direct consequence of the
previous definition
\cite[\S 1.2, \S 1.3, pp. 8–-17]{chi1}.
\begin{theorem} \label{theo2.3}
Let $\mathscr L$ be a moment functional and
$(p_n)$ a polynomial sequence. The
following statements are equivalent:
\begin{enumerate}
\item $(p_n)$is an OPS with respect to
$\mathscr L$.
\item ${\mathscr L}(\pi p_n)=0$ for all
polynomials $\pi$, $\deg \pi<n$, while
${\mathscr L}(\pi p_n)\ne 0$ if the $\deg \pi=n$.
\item ${\mathscr L}(x^mp_n(x))=K_n\delta_{n,m}$,
where $K_n\ne 0$, for $m=0, 1, \dots, n$.
\end{enumerate}
\end{theorem}
It is well-known that a monic OPS
$(p_n)$ satisfies a TTRR of the form
\eqref{1:2} where the coefficients
$\gamma_n$ do not vanish.
The converse is also true.
\begin{theorem} \label{theo2.4}(J. Favard)
Let $(p_n)$ be a polynomial sequence
satisfying the initial conditions
$p_{-1}=0$, $p_{0}=1$ and the TTRR
\eqref{1:2}, where $\gamma_n\neq 0$
for all $n\geq 1$.
Then $(p_n)$ is a OPS with respect
to the canonical moment functional
defined as
$$
\mathscr{L}(p_n)=\delta_{n,0},
\qquad n=0, 1, 2, \dots
$$
\end{theorem}
On the other side, if there exists
$\gamma_{N}=0$, then the sequence
$(p_n)$ can not be a OPS since the
identity
$$
\gamma_N=\frac{\mathscr{L}(p_N^2)}
{\mathscr{L}(p_{N-1}^2)},
$$
shows that condition (3) in definition
\ref{def2.2} does not hold.

Among the generalizations of OPS,
one is given by considering
symmetric bilinear functionals:
\begin{definition} \label{def2.5}
Given a sequence of polynomials $(p_n)$,
we say that $(p_n)$ is a OPS with
respect to a symmetric bilinear
functional $\mathscr B$ if the
following conditions hold:
\begin{enumerate}
\item $(p_n)$ is normal.
\item ${\mathscr B}(p_n,p_m)=0$, $m\ne n$.
\item ${\mathscr B}(p_n,p_n)\ne 0$.
\end{enumerate}
\end{definition}
It is usual to write $\langle f,g\rangle$
instead of $\mathscr{B}(f,g)$ when there
is no confusion about the bilinear
functional acting.
With this definition, the analog of
theorem \ref{theo2.3} is also valid but
it is not the  TTRR.

Notice that a sufficient condition for
the existence of the TTRR is the Hankel's
property, i.e.
$$
\langle t f,g\rangle=\langle f,tg\rangle,
$$
for all polynomials $f$ and $g$, where
$\langle\cdot,\cdot\rangle$ acts on
the variable $t$.

Among all the bilinear functionals we
focus on the following Sobolev type ones:
$$
\langle f,g\rangle=\mathscr{L}_0(fg)+
\mathscr{L}_1({\mathscr D}(f){\mathscr
D}(g)),
$$
where $\mathscr{L}_0$, $\mathscr{L}_1$
are linear functionals and ${\mathscr
D}$ is the derivative, the difference,
or the $q$-difference operator.

With all this overview we are going to
present an extension of Favard's theorem
for the case when some $\gamma$'s
coefficient vanishes.
%%%%%%%%%%%%%%%%%%%%%%%%%%%%%%%%%%%%%%%%%%%%%%%%%
\subsection{Degenerate generalization
of Favard's theorem} \label{sec2.1}
Consider the sequences $(\beta_n)$ and $(\gamma_n)$
of complex numbers and the polynomials generated
by the following recurrence relation:
\begin{equation} \label{aux}
p_{n+1}(x)=(x-\beta_n)p_n(x)-\gamma_n p_{n-1}(x),
\quad n=1, 2, \dots,
\end{equation}
with initial conditions $p_0(x)=1$ and
$p_1(x)=x-\beta_0$.
By Favard's result, we define for $n\ge 0$
the moment functional as
$$
{\mathscr L}_0(p_n)=\delta_{n,0}.
$$
Notice that in such a case, ${\mathscr L}_0(p_np_m)=0$
for all $n\ne m$.

It is important to point out
that if there exists $N$ so
that $\gamma_N=0$, then
$\mathscr{L}(p_N^2)=0$ and thus
the functional ${\mathscr L}_0$
does not determine the complete
polynomials sequence $(p_n)$.

In order to give an orthogonality
that characterizes the family
polynomials $(p_n)$, we need to
consider a linear operator
${\mathscr T}_1:\mathbb{P}
\longrightarrow \mathbb{P}$, and
polynomial sequence $(p_{n,1})$
satisfying the following conditions:
\begin{enumerate}
\item $\deg({\mathscr T}_1(p))=\deg(p)-1$
for any polynomial $p$.
\item The polynomial sequence
$(p_{n,1})$ is defined by
$$
p_{n-1,1}\stackrel{\rm  def}=
\frac{{\mathscr T}_1(p_n)}
{c_{n,1}},\qquad n\ge 1,
$$
where $c_{n,1}$ is the leading coefficient
of ${\mathscr T}_1(x^n)$,
and it satisfies, for $n\ge 1$, the
recurrence relation
$$
p_{n+1,1}(x)=(x-\beta_{n,1})p_{n,1}(x)
-\gamma_{n,1}p_{n-1,1}(x),
$$
where the sequence $(\gamma_{n,1})$ is such that
there exists a strictly increasing mapping
$$
\lambda: \{n\in \mathbb N:\gamma_{n,1}=0\} \longrightarrow
\{n\in \mathbb N:\gamma_n=0\}\,,
$$
with $\lambda(n)>n$.
\end{enumerate}
Observe that this last condition basically
means that, after the action of ${\mathscr
T}_1$, the possible vanishing $\gamma$'s
are shifted to a lower degree.
In fact, for many families of $q$-polynomials
and their relative natural $q$-difference
operator the condition about $\lambda$ writes
\begin{equation} \label{2:1}
\gamma_{n,1}=0 \iff  \gamma_{n+1}=0.
\end{equation}
Under these hypothesis, $(p_{n,1})$ is also
a monic orthogonal polynomials sequence with
respect to some moment functional, namely
${\mathscr L}_1$.
This procedure can be iterated $j$ times
giving a sequence of operators ${\mathscr
T}_j$, recurrence coefficients, $(\beta_{n,
j})$ and $(\gamma_{n,j})$, and moment
functionals ${\mathscr L}_j$ so that the
family
$$
({\mathscr T}_k\circ \cdots \circ {\mathscr
T}_2\circ {\mathscr T}_1(p_{n+k})),
$$
is orthogonal with respect to ${\mathscr
L}_{k}$.
If we denote by ${\mathscr T}^{(k)}
\stackrel{\rm def}={\mathscr T}_k
\circ \cdots \circ {\mathscr T}_1$, then
\begin{eqnarray}
\nonumber &p_{n,k}=C_{n,k}{\mathscr T}^{(k)}
(p_{n+k}),\qquad C_{n,k}\neq 0,\\ \label{2:2}
&p_{n+1,k}=(x-\beta_{n,k})p_{n,k}
-\gamma_{n,k}p_{n-1,k}, \ p_{0,k}=1, \ p_{1,k}=x-\beta_{0,k},\\
\nonumber  &{\mathscr L}_k(p_{n,k}p_{m,k})=0,
\quad n\neq m.
\end{eqnarray}
Taking into account this construction we
can state the degenerate generalization
of Favard's theorem.
\begin{theorem}\label{theo2.6}
Let $(p_n)$ be a polynomials sequence
satisfying the TTRR \eqref{aux},
so that there exists a
unique $N\in\mathbb{N}$ so
that $\gamma_N=0$, then $(p_n)$ is
the unique (monic) polynomial sequence
that fulfills the orthogonality
conditions
\begin{equation} \label{2:3}
\langle p_n,p_m\rangle=
{\mathscr L}_0(p_np_m)+
{\mathscr L}_{N}({\mathscr T}^{(N)}
(p_n){\mathscr T}^{(N)}(p_m))=0,
\qquad n\neq m.
\end{equation}
\end{theorem}
The choice of ${\mathscr T}^{(N)}$
and its link with ${\mathscr L}_{N}$
guarantees that $\langle p_n,p_m\rangle=0$
for all $n\ne m$. Hence, we only need to
check the orthogonality conditions $\langle
p_n,p_n\rangle \neq 0$ $\forall n\geq 0$ in
order to prove that $(p_n)$ is a MOPS (thus
the family $(p_n)$ is characterized by the
orthogonality property).

If $n<N$ then, by hypothesis,
$$
\langle p_n,p_n\rangle={\mathscr
L}_0(p_n^2)=\gamma_n\cdots\gamma_1\neq 0,
$$
and if $n\geq N$ then, taking into account \eqref{2:2},
$$
\begin{array}{rl}\langle p_n,p_n\rangle=& \displaystyle
{\mathscr L}_{N}({\mathscr T}^{(N)}(p_n){\mathscr T}^{(N)}
(p_n))=\frac{1}{C^2_{N,n-N}}{\mathscr L}_{N}(p_{n-N,N}
p_{n-N,N})\\& =\displaystyle \frac{\gamma_{n-N,N}
\gamma_{n-N-1,N}\dots\gamma_{1,N}}
{C^2_{n-N,N}}\neq 0. \qed\end{array}
$$

\begin{remark}
Notice that if there exists $N'<N$ such that
$\gamma_{n,N'}>0$, for all $n$, and the $\gamma$'s
coefficients satisfy $\gamma_1, \dots,\gamma_{N-1}>0$
then, the value $N$ in formula \eqref{2:3} can be
replaced by $N'$, and in such a case the proof of
the statement is similar.
Now $\langle p_n,p_n\rangle$ depends on the operators
$\mathscr{L}_0$ and $\mathscr{L}_{N'}$ and, in this case,
it is the sum of two positive terms which do not vanish
simultaneously.
\end{remark}

\begin{corollary} \label{cor2.7}
Let $(p_n)$ be a polynomial sequence
satisfying the TTRR \eqref{aux},
and let $\Lambda\stackrel{\rm def}=
\{n\in\mathbb{N}:\gamma_n=0\}$.
Then $(p_n)$ is the unique (monic) polynomial
sequence that fulfills the orthogonality
conditions
\begin{equation} \label{2:4}
\langle p_n,p_m\rangle={\mathscr L}_0(p_np_m)
+\sum_{k\in\mathscr{A}}{\mathscr L}_{k}({\mathscr
T}^{(k)}(p_n){\mathscr T}^{(k)}(p_m))=0,\qquad
n\neq m,
\end{equation}
being ${\mathscr A}=\{N_0,N_1,\dots\}$ with
$N_{j+1}=N_j+\min\{n:\gamma_{n,N_j}=0\}$.
\end{corollary}
The proof is straightforward taking
into account the proof of theorem \ref{theo2.6}.
\begin{remark}
Observe that if $\Lambda$ is a finite set,
then ${\mathscr A}$ is also finite as well.
Moreover, in \eqref{2:4} for any two polynomials
there is always a finite number of non
vanishing terms, so the result of corollary
\ref{cor2.7} remains even if the set $\mathscr A$ is
a infinite set.
\end{remark}

Among the operators $\mathscr T$ satisfying
the imposed conditions one of the most natural
ones is as follows:
$$
{\mathscr T}_1(p)(x)\stackrel{\rm def}=
{\mathscr L}_0\left(\frac{p(t)-p(x)}{t-x}
\right),\qquad\qquad p_{n,1}={\mathscr
T}_1(p_n)=p_n^{(1)},
$$
here ${\mathscr L}_0$ acts on the
variable $t$.
This is due the fact that this operator
commutes with the multiplication operator
by $x$, thus the coefficients of the TTRR
are shifted by 1, i.e., for $n\ge 1$,
$$
p_{n+1}^{(1)}(x)=(x-\beta_{n+1})p_n^{(1)}(x)
-\gamma_{n+1}p_{n-1}^{(1)}(x).
$$
Notice that, although theorem \ref{theo2.6}
seems to be new, it has been used implicitly
in \cite{cola1} where the operator
${\mathscr T}_j$ is the forward difference
operator $\Delta$ and it is applied to Racah,
Hahn dual Hahn, and Krawtchouk polynomials.
Orthogonality conditions of this type were
also used in \cite{AlAlRe,alpepire,mopepi,kwli1}
for Laguerre and Jacobi polynomials where
the operator $\mathscr T_j\equiv \mathscr T$
is the standard derivative, providing a
Sobolev type orthogonality to these families.
We are going to focus throughout this
paper on the orthogonality properties of
$q$-polynomials where the operator
$\mathscr T_j$ is a $q$-difference type
operator.
%%%%%%%%%%%%%%%%%%%%%%%%%%%%%%%%%%%%%%%%%%%%%%%%%
\subsection{The $q$-Askey and Nikiforov-Uvarov
tableaux} \label{sec2.2}
In this section we summarize the data for the
classical $q$-orthogonal families of the $q$-Hahn
tableau assuming that $\sigma(x)$ is a monic polynomial
(see, e.g., \cite{nisuuv,kost,mealma,koe1,alv5});
We also include the two families of $q$-polynomials
found by R. \'{A}lvarez-Nodarse and J. C. Medem
in \cite{alme}, namely the ``0-Jacobi/Bessel''
$q$-polynomials (0JB), and the ``0-Laguerre/Bessel''
$q$-polynomials (0LB) (see also cf.
\cite[pp. 214--217]{alv5}).
%%%%%%%%%%%%%%%%%%%%%%%%%%%%%%%%%%%%%%%%%%%%%%%%%
\begin{table}[!hbt]
\begin{tabular}{llll}
\hline
Family & Hyp. Repres. & Family & Hyp. Repres. \\
\hline
cd$q$H &${}_3\varphi_2(ae^{i\theta},
ae^{-i\theta};ab,ac|q)$ & b$q$J & ${}_3
\varphi_2(abq^{n+1},x;aq,cq|q)$ \\ $q$H &
${}_3\varphi_2(\alpha \beta q^{n+1},x;
\alpha q,q^{-N}|q)$ & d$q$H &${}_3\varphi_2
(q^{-x},\gamma \delta q^{x+1};\gamma q,q^{-N}|q)$
\\ 0JB & ${}_2\varphi_1(aq^n;-|x a^{-1}b^{-1})$ &
b$q$L & ${}_3\varphi_2(0,x;aq,bq|q)$ \\
l$q$J & ${}_2\varphi_1(abq^{n+1};aq|qx)$ &
$q$M & ${}_2\varphi_1(x;bq|-q^{n+1}c^{-1})$ \\
Q$q$K & ${}_2\varphi_1(x;q^{-N}|pq^{n+1})$ &
A$q$K & ${}_2\varphi_1(q^{-N}x;q^{-N}|xp^{-1})$\\
$q$K & ${}_2\varphi_1(x;xq^{N-n+1}|-pq^{n+N+1})$&
d$q$K & ${}_2\varphi_1(x;xq^{N-n+1}|cq x^{-1})$\\
0LB & ${}_2\varphi_1(0;-|x a^{-1})$ &
l$q$L & ${}_2\varphi_1(0;aq|qx)$ \\
$q$L & ${}_2\varphi_1(-x;0|q^{n+\alpha+1})$ &
A$q$C & ${}_2\varphi_1(-aq^n;0|qx)$ \\
$q$C  & ${}_2\varphi_1(x;0|-q^{n+1}a^{-1})$ &
ACI & ${}_2\varphi_1(x^{-1};0|qxa^{-1})$ \\
ACII & ${}_2\varphi_1(x;-|q^na^{-1})$ &
SW & ${}_1\varphi_1(;0|-q^{n+1}x)$ \\
\hline
\end{tabular}
\caption{\label{table2}Basic hypergeometric series
representation of some $q$-polynomials}
\end{table}
%%%%%%%%%%%%%%%%%%%%%%%%%%%%%%%%%%%%%%%%%%%%

Notice that the families for which $\sigma$
is not a polynomial on $x(s)$ are not included
in some of those tables; moreover to simplify
the notation, we use for table \ref{table2}
the following reduction:
$$
{}_r\varphi_s(\vec a;\vec b|z)\equiv {}_r
\varphi_s \left. \left(\begin{array}{c}q^{-n} \ \vec a\\
\vec b\end{array} \right|q;z \right).
$$

%%%%%%%%%%%%%%%%%%%%%%%%%%%%%%%%%%%%%%%%%%%%
\begin{table}[!hbt]
\begin{tabular}{llll}
\hline
$\sigma(x)$ & $(q-1)\tau(x)$ & $p_n(x)$ \\
\hline $(x-aq)(x-cq)$ & $(abq^2
-1)x+q(a+c-abq-acq)$ & $p_n(x;a,b,c;q)$ \\
$(x-1)(x-\alpha q^N)$ & $(\alpha \beta q^2-1)x+1-\alpha
(q-q^N+\beta q^{1+N})$ & $h_n^{(\alpha,\beta)} (x,N;q)$\\
$(x-aq)(x-bq)$ & $-x+q(a+b-abq)$ & $p_n(x;a,b;q)$ \\
$(x-q^{-N})(x-pq)$& $-x+q(q^{-N-1}+q-q^{-N+1})$
& $k_n^{Aff}(x;p,N;q)$\\
$(x-1)(x-a)$ & $-x+1+a$ & $u_n^{(\alpha)}(x;q)$
\\ $x(x-1)$ & $(abq^2
-1)x+1-aq$ & $p_n(x;a,b|q)$ \\$x(x-q^{-N})$ & $-(1+pq)x+(pq
-q^{-N})$& $k_n(x;p,N;q)$\\
$x(x-1)$ & $-(aq+1)x+1$ & $k_n(x;a;q)$  \\ $x(x-1)$ &
$-x+1-aq$ & $p_n(x;a|q)$ \\ $x^2$ & $(aq-1)x-abq$ &
$j_n(x;a,b)$ \\ $x^2$ & $-x+aq$ & $l_n(x;a)$ \\ $x-bq$ &
$qc^{-1}x-qc^{-1}-1+bq$ & $m_n(x;b,c;q)$ \\
$x-1$ & $-pq^{2-N}x+q{1-N}-1+pq$ & $k_n^{qtm}(x;p,N;q)$ \\
$x$ & $qx-1$ & $s_n(x;q)$ \\ $x$ & $\alpha qx-\alpha
q-1$& $l_n^{(\alpha)}(x;q)$ \\
$x$ & $qa^{-1}x-qa^{-1}-1$ & $c_n(x;a;q)$ \\
1 & $a^{-1}x-a^{-1}-1$ & $v_n^{(a)}(x;q)$ \\
\hline
\end{tabular}
\caption{\label{table1} Basic data of some monic
$q$-polynomials of the $q$-Hahn tableau}
\end{table}

%%%%%%%%%%%%%%%%%%%%%%%%%%%%%%%%%%%%%%%%%%%%%%%%%60%
\begin{figure}[!hbt]
\setlength{\unitlength}{0.00050000in}
\begingroup\makeatletter\ifx\SetFigFont\undefined%
\gdef\SetFigFont#1#2#3#4#5{%
  \reset@font\fontsize{#1}{#2pt}%
  \fontfamily{#3}\fontseries{#4}\fontshape{#5}%
  \selectfont}%
\fi\endgroup%
{\renewcommand{\dashlinestretch}{30}
\begin{picture}(9500,6750)(0,-10)
%%%legend
\drawline(519,6487)(1269,6487)
\dashline{60.000}(56,6237)(564,6242)
\drawline(464.251,6216.017)(564.000,6242.000)(463.759,6266.015)
\drawline(61,5975)(570,5975)
\drawline(470.099,5951.607)(570.000,5977.000)(469.903,6001.607)
%%%end legend
\put(556,6505){\makebox(0,0)[lb]{{\SetFigFont{7}{8.4}{\rmdefault}{\mddefault}{\updefault}Legend}}}
\put(766,6173){\makebox(0,0)[lb]{{\SetFigFont{7}{8.4}{\rmdefault}{\mddefault}{\updefault}particular case}}}
\put(739,5893){\makebox(0,0)[lb]{{\SetFigFont{7}{8.4}{\rmdefault}{\mddefault}{\updefault}limiting case}}}
%%% boxes
\drawline(4037,6651)(4987,6651)(4987,6164)
	(4037,6164)(4037,6651)
\drawline(2137,5226)(3087,5226)(3087,4739)
	(2137,4739)(2137,5226)
\drawline(5825,5226)(6980,5226)(6980,4739)
	(5825,4739)(5825,5226)
\drawline(3259,3801)(4290,3801)(4290,3314)
	(3259,3314)(3259,3801)
\drawline(4749,3801)(5750,3801)(5750,3314)
	(4749,3314)(4749,3801)
\drawline(7111,3801)(8061,3801)(8061,3314)
	(7111,3314)(7111,3801)
\drawline(949,3801)(1899,3801)(1899,3314)
	(949,3314)(949,3801)
\drawline(12,1913)(962,1913)(962,1426)
	(12,1426)(12,1913)
\drawline(1461,1913)(2411,1913)(2411,1426)
	(1461,1426)(1461,1913)
\drawline(2874,1913)(3824,1913)(3824,1426)
	(2874,1426)(2874,1913)
\drawline(4300,1913)(5250,1913)(5250,1426)
	(4300,1426)(4300,1913)
\drawline(5712,1913)(6662,1913)(6662,1426)
	(5712,1426)(5712,1913)
\drawline(8538,1913)(9488,1913)(9488,1426)
	(8538,1426)(8538,1913)
%%% iff lines
\dashline{50.000}(4290,3601)(4742,3601)
\drawline(4410,3621)(4290,3591)(4410,3561)
\drawline(4612,3631)(4742,3601)(4612,3561)
\dashline{50.000}(3830,1676)(4296,1676)
\drawline(3940,1706)(3820,1676)(3940,1646)
\drawline(4166,1706)(4296,1676)(4166,1646)
\drawline(4037,6651)(4987,6651)(4987,6164)
%%%% iff lines
%%% other lines
\dashline{60.000}(6375,4721)(3796,3811) %%%
\drawline(3876,3871)(3796,3811)(3926,3821)%%%
\drawline(4499,6160)(2616,5241)
\drawline(2694.903,5307.327)(2616.000,5241.000)(2716.833,5262.393)
\dashline{60.000}(4506,6160)(6374,5241)
\drawline(6273.235,5262.712)(6374.000,5241.000)(6295.307,5307.576)
\dashline{60.000}(2602,4734)(1435,3811)
\drawline(1497.925,3892.642)(1435.000,3811.000)(1528.942,3853.426)
\dashline{60.000}(2601,4729)(3783,3811)
\drawline(3688.687,3852.594)(3783.000,3811.000)(3719.356,3892.083)
\dashline{60.000}(2616,4721)(5188,3811)
\drawline(5085.388,3820.787)(5188.000,3811.000)(5102.065,3867.923)
\dashline{60.000}(2616,4728)(7555,3817)
\drawline(7452.124,3810.554)(7555.000,3817.000)(7461.194,3859.724)
\dashline{60.000}(3796,3317)(486,1918)
\drawline(568.378,1979.959)(486.000,1918.000)(587.843,1933.904)
\drawline(3796,3317)(1916,1925)
\drawline(1981.491,2004.598)(1916.000,1925.000)(2011.244,1964.414)
\drawline(5204,3307)(4770,1911)
\drawline(4775.814,2013.914)(4770.000,1911.000)(4823.560,1999.070)
\drawline(5215,3298)(6195,1931)
\drawline(6116.417,1997.707)(6195.000,1931.000)(6157.054,2026.839)
\dashline{60.000}(5215,3301)(7600,1925)
\drawline(7500.889,1953.319)(7600.000,1925.000)(7525.875,1996.628)
\dashline{60.000}(7574,3310)(1935,1928)
\drawline(2026.175,1976.085)(1935.000,1928.000)(2038.077,1927.522)
\dashline{60.000}(3796,3317)(3362,1921)
\drawline(3367.814,2023.914)(3362.000,1921.000)(3415.560,2009.070)
\drawline(7585,3302)(9017,1931)
\drawline(8927.479,1982.097)(9017.000,1931.000)(8962.056,2018.214)
\drawline(7580,3323)(6210,1921)
\drawline(6262.009,2009.995)(6210.000,1921.000)(6297.770,1975.050)
\drawline(7093,1912)(8043,1912)(8043,1425)
	(7093,1425)(7093,1912)
\drawline(6231,1418)(7557,503)
\drawline(7460.495,539.218)(7557.000,503.000)(7488.892,580.372)
\drawline(9058,1421)(7578,510)
\drawline(7650.055,583.709)(7578.000,510.000)(7676.265,541.129)
\drawline(7576,1411)(7578,531)
\drawline(7552.773,630.943)(7578.000,531.000)(7602.773,631.057)
\drawline(7101,499)(8051,499)(8051,12)
	(7101,12)(7101,499)
\put(4078,6303){\makebox(0,0)[lb]{{\SetFigFont{7}{8.4}{\rmdefault}{\mddefault}{\updefault}{AW$\equiv q$R}}}}
\put(2198,4890){\makebox(0,0)[lb]{{\SetFigFont{7}{8.4}{\rmdefault}{\mddefault}{\updefault}b$q$J$\equiv q$H}}}
\put(5890,4890){\makebox(0,0)[lb]{{\SetFigFont{7}{8.4}{\rmdefault}{\mddefault}{\updefault}cd$q$H$\equiv$d$q$H}}}
\put(1207,3473){\makebox(0,0)[lb]{{\SetFigFont{7}{8.4}{\rmdefault}{\mddefault}{\updefault}0JB}}}
\put(3301,3473){\makebox(0,0)[lb]{{\SetFigFont{7}{8.4}{\rmdefault}{\mddefault}{\updefault}b$q$L$\equiv$A$q$K}}}
\put(4792,3473){\makebox(0,0)[lb]{{\SetFigFont{7}{8.4}{\rmdefault}{\mddefault}{\updefault}$q$M$\equiv$Q$q$K}}}
\put(7228,3473){\makebox(0,0)[lb]{{\SetFigFont{7}{8.4}{\rmdefault}{\mddefault}{\updefault}l$q$J$\equiv$$q$K}}}
\put(287,1606){\makebox(0,0)[lb]{{\SetFigFont{7}{8.4}{\rmdefault}{\mddefault}{\updefault}0LB}}}
\put(1806,1606){\makebox(0,0)[lb]{{\SetFigFont{7}{8.4}{\rmdefault}{\mddefault}{\updefault}l$q$L}}}
\put(3176,1606){\makebox(0,0)[lb]{{\SetFigFont{7}{8.4}{\rmdefault}{\mddefault}{\updefault}ACI}}}
\put(4568,1606){\makebox(0,0)[lb]{{\SetFigFont{7}{8.4}{\rmdefault}{\mddefault}{\updefault}ACII}}}
\put(6040,1606){\makebox(0,0)[lb]{{\SetFigFont{7}{8.4}{\rmdefault}{\mddefault}{\updefault}$q$L}}}
\put(7415,1606){\makebox(0,0)[lb]{{\SetFigFont{7}{8.4}{\rmdefault}{\mddefault}{\updefault}$q$C}}}
\put(8783,1606){\makebox(0,0)[lb]{{\SetFigFont{7}{8.4}{\rmdefault}{\mddefault}{\updefault}A$q$C}}}
\put(7421,189){\makebox(0,0)[lb]{{\SetFigFont{7}{8.4}{\rmdefault}{\mddefault}{\updefault}SW}}}
\end{picture}
}\caption{Most of the known links between $q$-polynomials}
\end{figure}
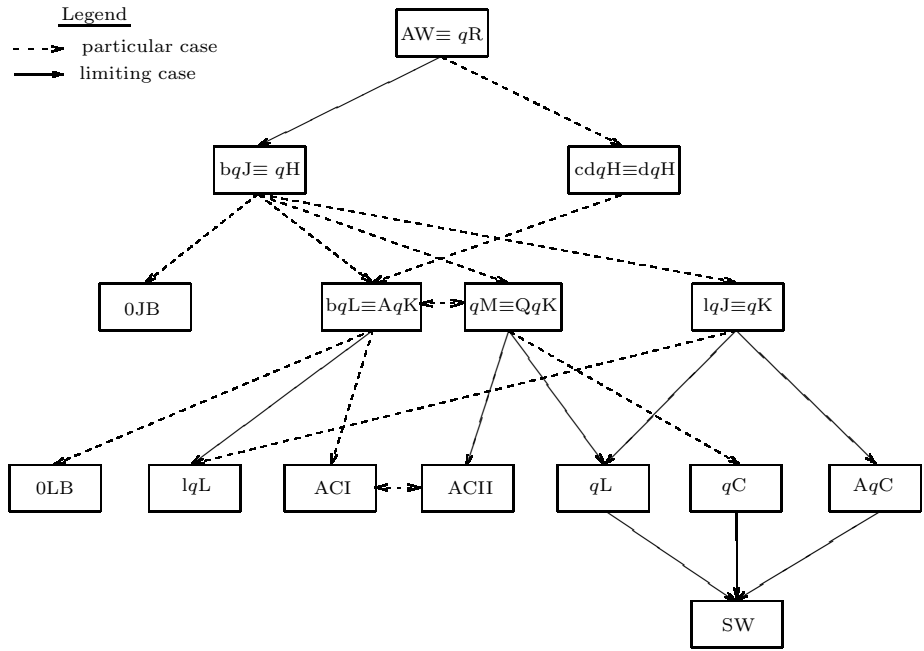
%%%%%%%%%%%%%%%%%%%%%%%%%%%%%%%%%%%%%%%%%%%%%%%%%60%

On the other hand, since the characterization
theorems characterize the $q$-po\-ly\-no\-mi\-als (see
e.g. \cite{alal,alv3,coma2}) then it is a direct
calculation by using table \ref{table1} to check
the following identities:
\begin{align*}
p_n(x+aq;a,a;q)=& l_n(x;a-a^2q),\\
p_n(x+aq;a,b,a;q)=&j_n(x;abq,1+ab^{-1}-b^{-1}
q^{-1}-aq).
\end{align*}
Although most of the identities we
present here are already known (see
\cite[\S 4]{kost}) we believe it
is a good idea to show them here
(see table \ref{table3}).
%%%%%%%%%%%%%%%%%%%%%%%%%%%%%%%%%%%%%%%%%%%%%%%%%
\section{The Askey-Wilson polynomials}\label{sec3}
This family of $q$-polynomials, which
were introduced by R. Askey and J.
Wilson in \cite{aswi1}, are located
at the top of the $q$-Askey tableau.
The monic Askey-Wilson polynomials
can be written as a basic hypergeometric
series \cite[p. 63]{kost}
$$
p_n(x;a,b,c,d|q)=
\frac{(ab;q)_n(ac;q)_n (ad;q)_n}
{(2a)^n (abcdq^{n-1};q)_n}{}_4\varphi_3
\left(\begin{array}{c|}q^{-n}, abcdq^{n-1}, ae^{i\theta},
ae^{-i\theta} \\ ab, ac, ad \end{array} \ q;q\right),
$$
with $x=\cos\theta$.
Moreover they fulfill, for $n\ge 0$, the TTRR
\begin{equation} \label{3:1}
xp_n(x)=p_{n+1}(x)+\beta_n p_n(x)+\gamma_n
p_{n-1}(x),
\end{equation}
where $\beta_n=\left(a+a^{-1}
-A_n-C_n\right)/2$, and $\gamma_n=
A_{n-1}C_n/4$ being
$$\begin{array}{rl}
A_n=& \displaystyle \frac{(1-abq^n)(1-acq^n)(1-adq^n)
(1-abcdq^{n-1})}{a(1-abcdq^{2n-1})
(1-abcdq^{2n})},\\ C_n=& \displaystyle
\frac{a(1-q^n)(1-bcq^{n-1})
(1-bdq^{n-1})(1-cdq^{n-1})}{(1-abcdq^{2n
-2})(1-abcdq^{2n-1})}.
\end{array} $$
Observe that, since $\Lambda=\{n\in\mathbb{N}:
\gamma_n=0\}$, then
$$
\Lambda=\emptyset \iff ab,ac,ad,bc,bd,cd\notin
\Omega(q)\stackrel{\rm def}=\{q^{-k}:k\in
\mathbb{N}_0\}.
$$
In the forthcoming sections
we only consider normal polynomials
sequences therefore $abcd \not \in
\Omega(q)$.
%%%%%%%%%%%%%%%%%%%%%%%%%%%%%%%%%%%%%%%%%%%%%%%%%
\subsection{The orthogonality conditions
for $|q|<1$} \label{sec3.1}

It is known that if the parameters
$a$, $b$, $c$, and $d$ are real, or
occur in complex conjugate pairs if
complex, $\max\{|a|,|b|,|c|,d|\}<1$,
the family fulfills the orthogonality
conditions
\cite{atsu4}
\begin{equation} \label{3:2}
\frac 1{2\pi}\int_{-1}^1 p_m(x)
p_n(x) \frac{\omega(x)}{\sqrt{1-x^2}}dx=
d_{n}^{2(AW)} \delta_{n,m},
\qquad n,m\ge 0,\end{equation}
where $d_{n}^{2(AW)}$ is the squared norm
of the monic Askey-Wilson polynomial of
degree $n$
\begin{equation}\label{3:3}
d_n^{2(AW)}=\frac{(abcdq^{2n};q)_\infty}
{4^n(abcdq^{n-1};q)_n(q^{n+1},abq^n,
acq^n,adq^n,bcq^n,bdq^n,cdq^n;q)_\infty},
\end{equation}
and
$$
\omega(x)=\left|\frac{(e^{2i\theta};q)_\infty}
{(ae^{i\theta},be^{i\theta},ce^{i\theta},d
e^{i\theta};q)_\infty}\right|^2=\frac
{h(x,1)h(x,-1)h(x,q^\frac 12)h(x,-q^\frac 12)}
{h(x,a)h(x,b)h(x,c)h(x,d)},
$$
with
$$
h(x,\alpha)\stackrel{\rm  def}=
\prod_{k=0}^\infty (1-2\alpha x
q^k+\alpha^2 q^{2k})=
(\alpha e^{i\theta},\alpha e^{-i
\theta};q)_\infty,\quad x=\cos\theta.
$$

Observe that the orthogonality conditions
given in \eqref{3:2} are a particular case
of the non-hermitian complex orthogonality
conditions
\begin{equation}\label{3:4}
\int_{\Gamma}p_n\left(\frac{z+z^{-1}}
{2}\right)p_m\left(\frac{z+z^{-1}}{2}
\right) W(z)dz=d_{n}^{2(AW)}
\delta_{n,m},\qquad n\neq m,
\end{equation}
which were obtained by Askey and Wilson
(see \cite{aswi1}), being
$$
W(z)=\frac{1}{z}w\left(\frac{z+z^{-1}}{2}\right).
$$
The poles of $w$ are
$$
\frac{\alpha q^k+(\alpha q^k)^{-1}}{2},
\quad \alpha=a,b,c,d,
\qquad k\in\mathbb{N}_0\,,
$$
therefore $W$ has convergent poles, since
$|q|<1$, at
$$
aq^k,\quad bq^k,\quad cq^k,\quad dq^k,\qquad
k\in\mathbb{N}_0,
$$
and divergent poles at
$$
a^{-1}q^{-k},\quad b^{-1}q^{-k},\quad
c^{-1}q^{-k},\quad d^{-1}q^{-k},\qquad
k\in\mathbb{N}_0.
$$
The contour $\Gamma$ is a curve
separating the divergent poles
from the convergent poles,
encircling them only once.
In fact, if the parameters satisfies
$\max\{|a|,|b|,|c|,|d|\}<1$ then
$\Gamma$ can be taken as the
unit circle, otherwise it is a
deformation of the unit circle.

The poles can be separated only if
$$
a^2,b^2,c^2,d^2,ab,ac,ad,bc,bd,cd
\notin \Omega(q)
,
$$
so in the following we focus our
attention when this does not occur.
Looking at the expression of the
coefficient $\gamma_n$, it
vanishes only if
$$
ab,ac,ad,bc,bd,cd\notin \Omega(q),
$$
and since any rearrangement of the
parameters does not change the
polynomial, it is enough to study
the following three key cases:
\begin{itemize}
\item
only $a^2\in\Omega(q)$, \item only
$ab\in\Omega(q)$, \item or only
$a^2=q^{-M}$ and $ab=q^{-N}$,
with $M<N-1$, belong to $\Omega(q)$.
\end{itemize}
%%%%%%%%%%%%%%%%%%%%%%%%%%%%%%%%%%%%%%%%%%%%%%%%%
\subsubsection{}\label{subsubsec311}
$a^2\in \Omega(q)$ and $ab,ac,ad,bc,bd,cd
\notin \Omega(q)$.
Although the poles can not be
separated, there is no $\gamma_n$
vanishing in the TTRR, so we look
for a simple reformulation of
\eqref{3:4}.
Let us assume $a^2=q^{-M}$ with
$M\in \mathbb N_0$, then the poles
that can not be separated are
$$
Z=\{q^{-M/2},\,  q^{1-M/2},  \dots,
\,q^{M/2}\}, \quad {\rm or}\quad
Z=\{-q^{-M/2},\, -q^{1-M/2},
\dots, \,-q^{M/2}\}.
$$
Notice that if some of these poles
coincide with the generated by $b$,
then $ab\in \Omega(q)$ which
is not possible in this case.
Hence $Z$ has empty intersection
with the rest of the poles of $W$.

We consider this case as the limit
for $p_n(\bullet;\alpha,b,c,d;q)$
with $\alpha\to a$, so the poles
of $W(\bullet;\alpha,b,c,d;q)$ can
be separated adequately.
Thus the orthogonality conditions
\eqref{3:4} are valid and can be
expressed as
$$
0=\int_{\Gamma'_1\cup\Gamma'_2}p_n
\left(\frac{z+z^{-1}}{2}\right)p_m
\left(\frac{z+z^{-1}}{2}\right)W(z)dz,
$$
where the curves $\Gamma'_1$ and
$\Gamma'_2$ separate the poles.
Therefore these curves can be
deformed in order to obtain the
integral through two curves,
$\Gamma_1$ and $\Gamma_2$, such
that they separate the convergent
poles from the divergent ones,
but the poles in $Z$ which
stand between the two curves,
with several residues added (see
next figures).

%%%%%%%%%%%%% figures 1 & 2
%%%%%%%%%%%%%%%%%%%%%%%%%%%%%%%% BEGIN FIGURE 1  %%%
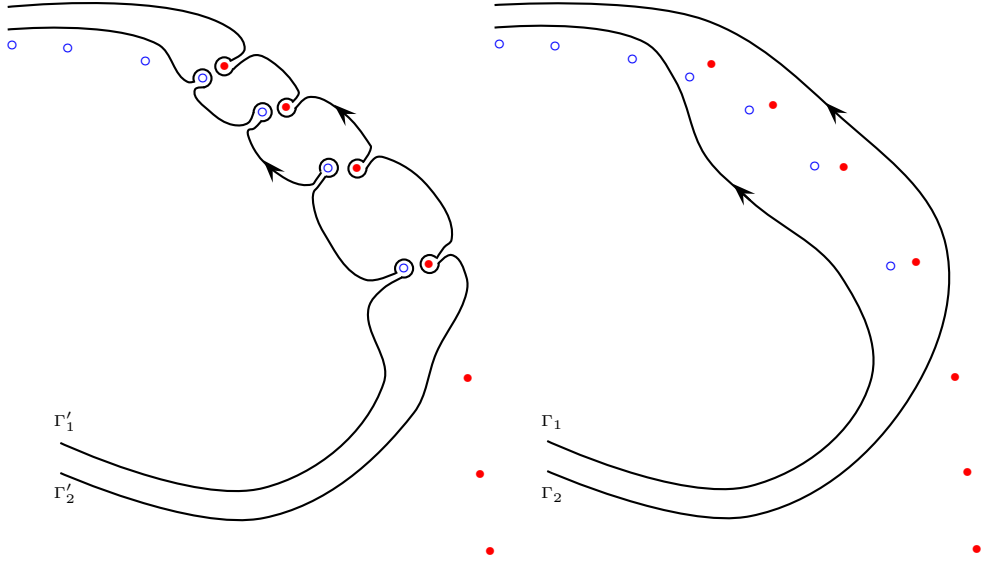
\setlength{\unitlength}{0.00045833in}

\begingroup\makeatletter\ifx\SetFigFont\undefined%
\gdef\SetFigFont#1#2#3#4#5{%
  \reset@font\fontsize{#1}{#2pt}%
  \fontfamily{#3}\fontseries{#4}\fontshape{#5}%
  \selectfont}%
\fi\endgroup%
{\renewcommand{\dashlinestretch}{30}
\begin{picture}(5569,6054)(300,200)
\thicklines
\thinlines
%%%%%%% red points (up 2 down)
\put(2480,5346){{\color{red}\circle*{90}}} % 1st
\put(3185,4874){{\color{red}\circle*{90}}} % 2nd
\put(3994,4171){{\color{red}\circle*{90}}} % 3rd
\put(4819,3077){{\color{red}\circle*{90}}} % 4th
\put(5262,1767){{\color{red}\circle*{90}}}
\put(5405,673){{\color{red}\circle*{90}}}
\put(5516,-208){{\color{red}\circle*{90}}}
%%%%%%% blue points
\put(53,5579){{\color{blue}\circle{77}}} % 1st
\put(690,5552){{\color{blue}\circle{77}}}  % 2nd
\put(1576,5400){{\color{blue}\circle{77}}} % 3rd
\put(2232,5204){{\color{blue}\circle{77}}} % 4th
\put(2914,4820){{\color{blue}\circle{77}}}
\put(3662,4179){{\color{blue}\circle{77}}}
\put(4530,3029){{\color{blue}\circle{77}}}
%%%%%%% arcs of circle (red)
\psarc(2.88,6.21){0.12}{45}{-13} % 1st
\psarc(3.71,5.66){0.12}{45}{5} % 2nd
\psarc(4.65,4.85){0.12}{45}{5} % 3rd
\psarc(5.61,3.58){0.12}{45}{5} % 4th
%%%%%%% arcs of circle (blue)
\psarc(2.59,6.05){0.12}{275}{205} % 1st
\psarc(3.39,5.6){0.12}{255}{185} % 2nd
\psarc(4.27,4.86){0.12}{250}{205} % 3rd
\psarc(5.27,3.52){0.12}{250}{205} % 4th
%  line below (C2) part 1
\pscurve(0.7,0.8)(3.4,0.2)(5.4,1.6)(5.7,2.4)(6.1,3.2)(6.1,3.4)(5.9,3.7)(5.72,3.58)
%  line below (C2) part 2
\pscurve(5.7,3.65)(5.8,3.8)(5.9,3.9)(5.6,4.6)(5,5)(4.9,5)(4.77,4.85)
%  line below (C2) part 3
\pscurve(4.72,4.93)(4.83,5.14)(4.82,5.16)(4.45,5.6)(4,5.8)(3.82,5.64)
\psline[arrowscale=2]{->}(4.45,5.6)(4.34,5.7)
%  line below (C2) part 4
\pscurve(3.78,5.73)(3.87,5.85)(3.86,5.89)(3.5,6.28)(3.3,6.36)(3.1,6.26)(3,6.17)
%  line below (C2) part 5
\pscurve(2.95,6.28)(3.15,6.44)(2.8,6.8)(0,7)
% line (C2) part 1
\pscurve(0.7,1.2)(3.4,0.6)(5,2)(4.8,3)(5,3.26)(5.24,3.4)
% line (C2) part 2
\pscurve(4.23,4.75)(4.12,4.6)(4.07,4.54)(4.07,4.34)(4.2,4)(4.76,3.37)(5.18,3.47)
% line (C2) part 3
\pscurve(3.37,5.48)(3.27,5.38)(3.2,5.36)(3.2,5.28)(3.5,4.85)(4.05,4.64)(4.1,4.7)(4.17,4.8)
\psline[arrowscale=2]{->}(3.5,4.85)(3.4,4.96)
% line (C2) part 4
\pscurve(2.62,5.93)(2.49,5.90)(2.7,5.6)(3,5.43)(3.2,5.46)(3.27,5.6)
% line (C2) part 5
\pscurve(2.49,6.01)(2.4,6)(2,6.5)(0,6.7)
%\psgrid(0,0)(0,0)(7,7)
\put(530,1170){\makebox(0,0)[lb]
{{\SetFigFont{6}{7.2}{\rmdefault}
{\mddefault}{\updefault}$\Gamma'_1$}}}
\put(530,341){\makebox(0,0)[lb]
{{\SetFigFont{6}{7.2}{\rmdefault}
{\mddefault}{\updefault}$\Gamma'_2$}}}
\end{picture}
}
%%%%%%%%%%%%%%%%%%%%%%%%%%%%%%%%%%%%%%%%%%%%%%%%%%%%%%
\setlength{\unitlength}{0.00045833in}

\begingroup\makeatletter\ifx\SetFigFont\undefined%
\gdef\SetFigFont#1#2#3#4#5{%
  \reset@font\fontsize{#1}{#2pt}%
  \fontfamily{#3}\fontseries{#4}\fontshape{#5}%
  \selectfont}%
\fi\endgroup%
{\renewcommand{\dashlinestretch}{30}
\begin{picture}(5569,6054)(-5269,-5904)
\thicklines
%  line below (C2)
\pscurve(0.7,0.8)(3.4,0.2)(5.4,1.6)(6,3.8)(4.5,5.6)(2.8,6.8)(0,7)
\psline[arrowscale=2]{->}(4.5,5.6)(4.4,5.7)
% line (C1)
\pscurve(0.7,1.2)(3.4,0.6)(5,2)(4.6,3.4)(3.4,4.4)(2.8,5)(2.4,6)(2,6.5)(0,6.7)
\psline[arrowscale=2]{->}(3.4,4.4)(3.17,4.62)
\thinlines
\put(3185,4874){{\color{red}\circle*{90}}}
\put(2232,5204){{\color{blue}\circle{90}}}
\put(2480,5346){{\color{red}\circle*{90}}}
\put(4819,3077){{\color{red}\circle*{90}}}
\put(2914,4820){{\color{blue}\circle{90}}}
\put(5516,-208){{\color{red}\circle*{90}}}
\put(53,5579){{\color{blue}\circle{90}}}
\put(690,5552){{\color{blue}\circle{90}}}
\put(4530,3029){{\color{blue}\circle{90}}}
\put(3994,4171){{\color{red}\circle*{90}}}
\put(5405,673){{\color{red}\circle*{90}}}
\put(3662,4179){{\color{blue}\circle{90}}}
\put(1576,5400){{\color{blue}\circle{90}}}
\put(5262,1767){{\color{red}\circle*{90}}}
\put(530,1170){\makebox(0,0)[lb]
{{\SetFigFont{6}{7.2}{\rmdefault}
{\mddefault}{\updefault}$\Gamma_1$}}}
\put(530,341){\makebox(0,0)[lb]
{{\SetFigFont{6}{7.2}{\rmdefault}
{\mddefault}{\updefault}$\Gamma_2$}}}
\end{picture}
}
%%%%%%%%%%%%%%%%%%%%%%%%%%%%%%% END FIGURE 2  %%%
%%%%%%%%%%%%% end figures 1 & 2
\vspace{-60mm}

When $\alpha\to a$, the poles $\alpha
q^k$ with $k\leq M$ and $\alpha^{-1}
q^{-(M-k)}$ converges to $aq^{k}$ and
it can be seen that the sum of the
two residues at this points tends to zero.
So the limit $\alpha\to a$ yields
$$\int_{\Gamma_1\cup \Gamma_2}p_n
\left(\frac{z+z^{-1}}{2}\right)
p_m\left(\frac{z+z^{-1}}{2}\right)
W(z)dz=d_n^2\delta_{n,m},
$$
with $d_n^2$ the normalizing factor
given by \eqref{3:3}.

\subsubsection{}\label{subsubsec312}
$ab=q^{-N+1}$ and $a^2, b^2\notin\{q^0,
\dots,q^{-N+2}\}$, i.e $\gamma_N=0$,
so $N\in \Lambda$.
The orthogonality conditions depend
on the size of $\Lambda$ (see corollary
\ref{cor2.7}), so we show how it is in the
simplest case $\Lambda=\{N\}$, i.e.,
$ac,ad,bc,bd,cd\notin\Omega(q)
\setminus\{q^{-N}\}$.

Since monic $q$-Racah polynomials can be
written in terms of the basic hypergeometric
functions as \cite[(3.2.1)]{kost}
$$
r_n(\mu(x);\alpha,\beta,\gamma,\delta|q)=
{}_4\varphi_3\left.\left(\begin{array}{c}
q^{-n}, \alpha\beta q^{n+1},q^{-x},\gamma\delta
q^{x+1}\\ \alpha q,\beta\delta q,\gamma q
\end{array} \right|q;q\right),
$$
with $\mu(x)=q^{-x}+\gamma\delta q^{x+1}$,
the following identity linking
Askey-Wilson and $q$-Racah polynomials
holds
$$
p_n(x;a,b,c,d;q)=r_n(2ax;q^{-N},
cdq^{-1},adq^{-1},ad^{-1};q),
$$
and it yields the moment
functional $\mathscr{L}_0$ in
theorem \ref{theo2.6} which is the
one known for $q$-Racah polynomials
\begin{equation}\label{3:5}
\mathscr{L}_0(p)=\sum_{j=0}^{N-1}
\frac{(q^{-N+1},ac,ad,a^2;q)_j}
{(q,a^2q^N,ac^{-1}q,ad^{-1}q;q)_j}
\frac{(1-a^2q^{2j})}
{(cdq^{-N})^j(1-a^2)} p\left(
\frac{q^{-j}+a^2q^j}{2a}\right)\,.
\end{equation}
Notice that the assumptions on $a^2$
and $b^2$ guarantees the
definition of $\mathscr{L}_0$.

Furthermore, since
$$
\mathscr{D}_q p_n(x;a,b,c,d;q)=
\frac{q^n-1}{q-1} p_{n-1}(x;aq^{1/2},
bq^{1/2},cq^{1/2},dq^{1/2};q),
$$
where the $q$-difference operator,
also called the Hahn's operator, is
$$
{\mathscr D}_q(f)(z)
\stackrel{\rm def}=\left\{
\begin{array}{l}\displaystyle \frac{f(z)-f(qz)}
{(1-q)z},\quad z\ne 0 \ \wedge
\ q\ne 1,\\
f'(z),\quad z=0\ \vee \ q=1,\end{array}
\right.
$$
the operator $\mathscr{T}$ can
be chosen as $\mathscr{D}_q$ and
the condition \eqref{2:1} holds.
Hence, for $n\geq N$,
$$
\mathscr{D}_q^N p_n(x;a,b,c,d;q)=
\frac{(q^{n-N+1};q)_N}{(1-q)^N}
p_{n-N}(x;aq^{N/2},bq^{N/2},
cq^{N/2},dq^{N/2};q),
$$
so $\mathscr{L}_N$ is the moment
functional associated with the
Askey-Wilson polynomials with
parameters
$aq^{N/2}$, $bq^{N/2}$, $cq^{N/2}$,
and $dq^{N/2}$, i.e.
$$
\mathscr{L}_N(p)=\int_{\Gamma} p
\left(\frac{z+z^{-1}}{2}\right)\,
\frac{1}{z}w\left(\frac{z+z^{-1}}
{2}\right)\,dz,
$$
where
$$
w(z)=w(z;aq^{N/2},bq^{N/2},cq^{N/2},dq^{N/2};q),
$$
and $\Gamma$ is a contour which separates
the poles.
Then, by theorem \ref{theo2.6}, the polynomial
sequence $(p_n(x;a,b,c,d))$ is uniquely
determined, up to a constant, by the orthogonality
conditions, for $n\ne m$,
$$
\langle p_n(\bullet;a,b,c,d;q), p_m(\bullet;a,b,
c,d;q)\rangle=
\mathscr{L}_0 (p_np_m)+\mathscr{L}_N({\mathscr
D}_{q}^{N}(p_n){\mathscr D}_{q}^{N}(p_m))=0.$$

\subsubsection{} \label{subsubsec313}
$ab=q^{-N+1}$ and $a^2=q^{-M}$, with $M\in\{0,
\dots,N-2\}$.
Also the form of the orthogonality depends on
the numbers of elements of $\Lambda$.
For simplicity, we see only the case when the
cardinal of $\Lambda$ is one, and when $\Lambda$
is greater, the orthogonality is given by corollary
\ref{cor2.7}.

The orthogonality is basically the same that in
the case \ref{subsubsec312}, but now $\mathscr{L}_0$
is not valid since it has lost several orthogonality
conditions.
The adequate form of $\mathscr{L}_0$ is obtained
as a limit case.
Let us consider the linear functional
$$
\mathscr{L}_0^\alpha(p)=\sum_{j=0}^{N-1}A_j(\alpha)
p\left(\mu_j(\alpha)\right),
$$
with $\mu_j(x;\alpha)=(\alpha q^{j}+\alpha^{-1}q^{-j})/2$, and
$$
A_j(\alpha)=\frac{(q^{-N+1},\alpha c,\alpha d,
\alpha^2;q)_j} {(q,\alpha^2q^N,\alpha c^{-1}q,\alpha
d^{-1}q;q)_j}\, \frac{(1-\alpha^2q^{2j})}{(cdq^{-N})^j
(1-\alpha^2)}.
$$
Straightforward computations yields
$$
A_j(a)=0,
$$
for $j\in\{M+1,\dots,N-1\}$ and $j=M/2$
if $M$ is even, and
$$
A_j(a)+A_{M-j}(a)=0,\qquad
\mu_{j}(a)=\mu_{M-j}(a),
$$
for $j\in\{0,\dots, M\}$ but $j=M/2$
if $M$ is even.
Thus $\mathscr{L}_0^\alpha$ tends
to the null functional.
But since it is possible to
consider any normalization,
we remove the common factor
$(\alpha-a)$,
$$
\lim_{\alpha\to a} \frac{A_j(\alpha)}
{\alpha-a} p(\mu_j(\alpha))=A_j'(a)
p(\mu(a)),
$$
for $j=M+1,\dots,N$ and if $M$ is even
$j=M/2$, and also
\begin{align*}
&\lim_{\alpha\to a}\frac{A_j(\alpha)
p(\mu_j(\alpha))+A_{M-j}(\alpha)
p(\mu_{M-j}(\alpha))}{\alpha-a}\\&
\qquad=(A_j'(a)+A_{M-j}'(a))
p(\mu_j(a))+A_j(a)(q^j-q^{M-j})
p'(\mu_j(a))
\end{align*}
for $j=0,\dots, M$ but if $M$ is even
$j\neq M/2$.

Hence we define $\mathscr{L}_0$ as
\begin{align*}
\mathcal{L}_0(p)&=\sum_{j=0}^{(M-1)/2}
(A_j'(a)+A_{M-j}'(a)) p(\mu_j(a))+A_j(a)
(q^j-q^{M-j})p'(\mu_j(a))\\&\quad+
\sum_{j=M+1}^{N-1} A_j'(a) p(\mu_j(a))
\end{align*}
 if $M$ is odd, and
\begin{align*}
\mathcal{L}_0(p)&=\sum_{j=0}^{M/2-1}
(A_j'(a)+A_{M-j}'(a)) p(\mu_j(a))+
A_j(a)(q^j-q^{M-j})p'(\mu_j(a))\\&
\quad+\sum_{j=M+1}^{N-1} A_j'(a)
p(\mu_j(a))+A_{M/2}'(a) p(\mu_{M/2}(a))
\end{align*}
if $M$ is even.
The Askey-Wilson polynomials
of degree at most $N$ with $ab=q^{-N+1}$,
$a^2=q^{-M}$ and $M=0,\dots N-2$ are uniquely
determined by the orthogonality property:
$$
\mathscr{L}_0(p_np_m)=0,\qquad 0\le m<n
\leq N.
$$
In particular, $p_N$ has simple roots
on $\mu_j(a)$, $j=M+1,\dots, N$ and
on $\mu_{M/2}(a)$ if $M$ is even; the
rest of the roots, $\mu_{j}(a)$, $j=0,
\dots,[(M-1)/2]$ are double.

The moment functional $\mathscr{L}_N$ in
theorem \ref{theo2.6} is the same that the
one given in section \ref{subsubsec312}.

%%%%%%%%%%%%%%%%%%%%%%%%%%%%%%%%%%%%%%%%%%%%%%%%%%
\subsection{The orthogonality conditions
for $|q|\geq 1$} \label{sec3.2} Taking into account the
relation between basic hypergeometric series \cite[p. 9]{kost}
\begin{equation} \label{3:6}
{}_4\varphi_3\left.\left(\begin{array}{c} q^{-n},a,b,c\\d,e,f
\end{array} \right|q;q\right)={}_4\varphi_3\left.\left(\begin{array}{c}
q^{n},a^{-1},b^{-1},c^{-1}\\d^{-1},e^{-1},f^{-1}
\end{array} \right|q^{-1};\frac{abcq^{-n}}{def}\right).
\end{equation}
We can relate each family of $q$-polynomials
on the parameter $q$ into another family
of $q$-polynomials on the parameter $q^{-1}$. In fact in this case it provides
$$
p_n(x;a,b,c,d|q^{-1})=p_n(x;a^{-1},b^{-1},
c^{-1},d^{-1}|q).
$$
Therefore if $|q|>1$ we can get analogous
orthogonality conditions just using this
relation and the orthogonality conditions
given in Section \ref{sec3.1} for $|q|<1$.

If $q$ is a primitive root of unity, i.e.
$q=e^{2\pi i M/N}$ with $\gcd(N,M)=1$ then
$\{k N:k\in\mathbb{N}\}\subseteq\Lambda$,
so, by corollary \ref{cor2.7}, for $k=1$
we need to construct the following orthogonality
property for the Askey-Wilson polynomials
up to degree $N$ \cite{atk}, i.e.
\begin{equation} \label{3:7}
\sum_{s=0}^{N-1} p_n(x_s)p_m(x_s)\omega_s=
\gamma_1\cdots \gamma_n\delta_{n,m},
\end{equation}
where $n, m=0, 1, \dots ,N-1$, $\{x_s\}_{s=0}^{N-1}$
are the zeroes of $p_N$, and the weight function is
$$
\omega_s=\frac{\gamma_1\cdots \gamma_{N-1}}
{p_{N-1}(x_s)p'_N(x_s)}.
$$
Observe that the only requirement to be added
is that all zeros $x_s$ must be simple.

Since the method considered in \cite{Zhedanov}
to obtain $\omega_s$ can be applied to obtain
such weights functions to other families of
$q$-polynomials, next we give a brief outline of it.

It is known that Askey-Wilson polynomials are
polynomial eigenfunctions of the second order
homogeneous linear difference operator:
$$
\sigma(-s)\frac{\Delta p_n(x(s))}{\Delta x(s)}+
\sigma(s)\frac{\nabla p_n(x(s))}{\nabla x(s)}-
\lambda_n\Delta x(s-\mbox{$\frac 12$})p_n(x(s))=0,
$$
being $\sigma(s)=-(q^{1/2}-q^{-1/2})^2
q^{-2s+1/2}(q^s-a)(q^s-b)(q^s-c)(q^s-d)$, and their
corresponding eigenvalues
$$
\lambda_n=-4q^{-n+1}(1-q^n)(1-abcdq^{n-1}).
$$

Notice that such difference operator can be
rewritten, by using the definition of the
difference operators $\Delta$ and $\nabla$,
as \cite[Eq. (3.1.7)]{kost}, \cite[Eq. (3.6)]{Zhedanov}
$$
A(z^{-1}) p_n(q^{-1}z)-(A(z)+A(z^{-1})) p_n(z)+
A(z) p_n(qz)=\lambda_n p_n(z),
$$
where $A(z)=(1-az)(1-bz)(1-cz)(1-dz)/((1-z^2)(1-qz^2))$.
Therefore, multiplying the previous equation by a
function $\rho(s)$ satisfying the only requirement
of periodicity $\rho(s+N)=\rho(s)$, and combining
it with a similar equation for the polynomials
$p_m(x_s)$, one can get a bilinear relation:
$$\begin{array}{c}
A_s\sigma(s)\left(p_n(x_{s-1})p_m(x_s)-p_n(x_{s})
p_m(x_{s-1})\right)\\ +C_s\sigma(s)\left(p_n(x_{s+1})
p_m(x_s)-p_n(x_{s})p_m(x_{s+1})\right)\\
=(\lambda_n-\lambda_m)\sigma(s)p_n(x_s)p_m(x_s).
\end{array}$$

Choose $\rho(s)$ in such a way that
\begin{equation} \label{3:8}
A_{s+1}\rho(s+1)=C_s\rho(s),
\end{equation}
summing from $s=0$ to $s=N-1$ and using the
obvious periodicity property of $\rho(s)$ we get
the orthogonality property:
$$
(\lambda_n-\lambda_m)\sum_{s=0}^{N-1} p_n(x_s)
p_m(x_s)\rho(s)=0, \qquad n\ne m.
$$
Hence $\omega_s=\omega_0\rho(s)$, with $\omega_0$
is the normalization constant, is determined from
the relation \eqref{3:8}.

Spiridonov and Zhedanov found that the
polynomials $(p_n(\bullet;a,b,c,d;e^{2\pi i
M/N}))$, with $0\le n\le N$, under the assumptions
$$
abcd,ab,ac,ad,bc,bd,cd\neq q^{k},
\qquad k=0,\dots, N-1\,,
$$
are uniquely determined by the
orthogonality conditions
$$
\mathscr{L}_0(p_np_m)=d_n^2 \delta_{n,m},
\quad d_n^2\ne 0,
$$
being
$$
\mathscr{L}_0(p)=\sum_{j=0}^{N-1}\left(
\frac{q}{abcd}\right)^j\frac{(1-rq^{2j})
(ar,br,cr,dr;q)_j}{(1-r^2)(qr/a,qr/b,qr/c,
qr/d;q)_j}p\left(rq^j+r^{-1}q^{-j}\right),
$$
and $r$ the root with minimal argument
of the equation
$$
r^N=E_N/2+\sqrt{E_N^2/4-1},
$$
being
$$
E_N=\frac{a^N+b^N+c^N+d^N-(abc)^N-(abd)^N-
(acd)^N-(bcd)^N}{1-(abcd)^N}.
$$
\begin{remark}  \label{rem5.2}
A straightforward computation shows that
$$
\rho(s)\stackrel{\rm def}=\left(
\frac{q}{abcd}\right)^s\frac{(1-rq^{2s})
(ar,br,cr,dr;q)_s}{(1-r^2)(qr/a,qr/b,qr/c,
qr/d;q)_s},
$$
satisfies the condition \eqref{3:8}.
A hint for such calculation can be
found in \cite[Lemma 5.1]{coma2}.
\end{remark}

Due to the cyclic behavior of the TTRR coefficients and
since $\gamma_N=0$, these polynomials satisfy the
identity
$$
p_n=p_N^\ell p_m,\qquad n=\ell N+m,\qquad 0\leq m<N,
$$
which explains the behavior of the polynomial
for greater degrees.
However corollary \ref{cor2.7} is applicable.
For $n\geq N$
$$%\begin{array}{rl}
\mathscr{D}^{N}_q p_n(x;a,b,c,d;q)=
%&\displaystyle \frac{(q^{n-N+1};q)_N}{(1-q)^N}
%p_{n-N}(x;aq^{N/2},bq^{N/2},cq^{N/2},dq^{N/2};q)\\=&
\displaystyle \frac{(q^{n-N+1};q)_N}{(1-q)^N}
p_{n-N}((-1)^Mx;a,b,c,d;q), %\end{array}
$$
so the orthogonality conditions that characterizes
all polynomials are the following:
\begin{itemize}
\item If $M$ is even:
$$
\langle p_n,p_m\rangle=\sum_{j=0}^{\infty}
\mathscr{L}_0(\mathscr{D}_q^{Nj}(p_n)
\mathscr{D}_q^{Nj}(p_m)).
$$l
\item If $M$ is odd:
$$
\langle p_n,p_m\rangle=
\sum_{j=0}^{\infty}\mathscr{L}_0(\mathscr{D}_q^{2jN}(p_n)
\mathscr{D}_q^{2jN}(p_m))+\mathscr{L}_N(
\mathscr{D}_q^{(2j+1)N}(p_n)\mathscr{D}_q^{(2j+1)N}(p_m)),
$$
being
$$
\mathscr{L}_N(p)=\sum_{j=0}^{N-1}\left(
\frac{q}{abcd}\right)^j\frac{(1-rq^{2j})
(ar,br,cr,dr;q)_j}{(1-r^2)(qr/a,qr/b,qr/c,
qr/d;q)_j}p\left(-rq^j-r^{-1}q^{-j}\right).
$$
\end{itemize}
%%%%%%%%%%%%%%%%%%%%%%%%%%%%%%%%%%%%%%%%%%%%%%%%
\section{The big $q$-Jacobi polynomials} \label{sec4}
The big $q$-Jacobi polynomials, which were
introduced by Hahn in 1949, are located at the
top of the $q$-Hahn tableau.
The monic big $q$-Jacobi polynomials can be
written in terms of basic hypergeometric series
as \cite[p. 73]{kost}
\begin{equation} \label{4:1}
p_n(x;a,b,c;q)=\frac{(aq,cq;q)_n}
{(abq^{n+1};q)_n}{}_3\varphi_2\left.\left(\begin{array}{c}
q^{-n},abq^{n+1},x\\ aq,cq\end{array}\right|q; q\right).
\end{equation}

In fact they are the most general family of
$q$-polynomials on the $q$-exponential
lattice, also called $q$-linear lattice;
and they appear, among others branches of
physics, in the representation theory of the
quantum algebras \cite{vikl}.
The monic big $q$-Jacobi polynomials fulfill,
for $n\ge 1$, the following TTRR:
\begin{equation} \label{4:2}
xp_n(x)=p_{n+1}(x)+\beta_n p_n(x) +\gamma_n
p_{n-1}(x),
\end{equation}
with $\beta_n=1-\hat A_n-\hat C_n$, and
$\gamma_n=\hat A_{n-1} \hat C_n$ being
\begin{equation}\label{4:3}\begin{array}{rl}
\hat A_n=& \displaystyle \frac{(1-aq^{n+1}) (1-abq^{n+1})
(1-cq^{n+1})}{(1-abq^{2n+1})(1-abq^{2n+2})},\\
\hat C_n=& \displaystyle -acq^{n+1}\frac{(1-q^n)
(1-abc^{-1}q^n)(1-bq^n)}{(1-ab q^{2n})
(1-abq^{2n+1})}.
\end{array}
\end{equation}

\begin{remark}  \label{rem4.1}
Observe that if $a=0$ then the coefficients $\gamma_n=0$
for all $n\in \mathbb N_0$, and if $b=0$, or $c=0$, then
the big $q$-Jacobi polynomials become the big
$q$-Laguerre or the little $q$-Jacobi polynomials
respectively, which are located below in the $q$-Askey
tableau thus we omit these cases.
\end{remark}

A slightly less detailed study on orthogonality
conditions for the big $q$-Jacobi can be found in
\cite{moga}.
%%%%%%%%%%%%%%%%%%%%%%%%%%%%%%%%%%%%%%%%%%%%%%%%%
\subsection{The orthogonality conditions
for $|q|<1$} \label{sec4.1}

It is known that if $0<q<1$, $0<a,b<q^{-1}$, and
$c<0$ the family of big $q$-Jacobi polynomials
fulfills the orthogonality conditions
\cite[p. 73]{kost}
\begin{equation} \label{4:4}
\int_{cq}^{aq} \frac{(a^{-1}x,c^{-1}x;q)_\infty}
{(x,bc^{-1}x;q)_\infty}p_n(x;a,b,c;q)
p_m(x;a,b,c;q)d_qx= d_n^{2(BqJ)}\delta_{n,m},
\end{equation}
where the Jackson $q$-integral (see \cite{gara,kost})
is defined as follows
$$
\int_{a}^bf(t)d_qt=a(q-1)\sum_{s=0}^\infty
f(aq^{s})q^s-b(q-1)\sum_{s=0}^\infty f(bq^{s})q^s.
$$
The aim of this section is to give orthogonality
conditions for the big $q$-Jacobi polynomials for
general complex parameters, including complex $|q|<1$,
except for those for which the family is not normal,
i.e. $ab\in\Omega(q)$.

In fact, notice that if the parameters belong to
compact sets where the integrand in \eqref{4:4}
is bounded, hence such series converges uniformly.
Thus we can apply the Weierstrass theorem and
analytic prolongation in order to asserts that
\eqref{4:4} is valid for
$$
a,b,c,abc^{-1}\notin  \Omega(q),
$$
which it is equivalent to
$\Lambda=\emptyset$, therefore in the
following we focus our attention in the case
$\Lambda\not= \emptyset$.
More precisely, we study the cases for which
this set has exactly one element, namely $N$.
If this set is greater we refer the reader to
corollary \ref{cor2.7}.

\subsubsection{}\label{subsub411}
$c=q^{-N}$  and $a,b,abc^{-1}\notin \Omega(q)\setminus\{q^{-N}\}$.
Taking into account that the big $q$-Jacobi and $q$-Hahn
polynomials are linked through the relation
$$
p_n(x;a,b,q^{-N};q)=h_n^{(a,b)}(x;N-1;q),
$$
the moment functional ${\mathscr L}_0$ in theorem
\ref{theo2.6} is the one known for the $q$-Hahn
polynomials \cite{kost} with parameters $a$, $b$
and $N-1$.
\begin{equation}\label{4:5}
{\mathscr L}_0(p)=
\sum_{x=0}^{N-1}\frac{(aq,q^{-N+1};q)_x}
{(q,b^{-1}q^{-N+1};q)_x}(abq)^{-x}~p(q^{-x}).
\end{equation}
Moreover, since
$$
{\mathscr D}_{q^{-1}}  h_n^{(\alpha,\beta)}(x;M;q)=
\frac{q^{-n}-1}{q^{-1}-1} h_{n-1}^{(\alpha q,\beta q)}(x;M-1;q),
$$
the operator $\mathscr{T}$ in theorem
\ref{theo2.6} can be chosen as
${\mathscr D}_{q^{-1}}$ and the condition
\eqref{2:1} holds (see the relation
between $q$-Hahn and big $q$-Jacobi polynomials
and the expression \eqref{4:3} for the
coefficients $\gamma_n$).
Also, for $n\geq N$,
$$
{\mathscr D}_{q^{-1}}^{N}p_n(x;a,b,q^{-N};q)=
\frac{(q^{-n};q)_N}
{(1-q^{-1})^N}p_{n-N}(x;aq^N,bq^N,1;q).
$$
Accordingly with these expressions and the weight
function for the big $q$-Jacobi polynomials with
parameters $aq^N$, $bq^N$ and $1$, if we define
$$
\mathscr{L}_N(p)=\int_q^{aq^{N+1}}\frac{(a^{-1}
q^{-N}x;q)_\infty}{(bq^Nx;q)_\infty} p(x)~d_qx,
$$
then, by theorem \ref{theo2.6}, the orthogonality
conditions for $n\ne m$,
$$
\langle p_n(\bullet;a,b,q^{-N};q), p_m(\bullet;a,b,
q^{-N};q)\rangle=
\mathscr{L}_0 (p_np_m)+\mathscr{L}_N({\mathscr
D}_{q^{-1}}^{N}(p_n){\mathscr D}_{q^{-1}}^{N}
(p_m))=0,
$$
determine uniquely the big $q$-Jacobi polynomials
for all non-negative integer degrees up to a
constant factor.

\subsubsection{} \label{subsub412} $a=q^{-N}$,
$b,c,abc^{-1}\notin\Omega(q)\setminus\{q^{-N}\}$.
By using the identity
\begin{equation} \label{4:6}
p_n(x;a,b,c;q)=p_n(x;c,abc^{-1},a;q),
\end{equation}
which can be obtained easily from the
hypergeometric representation \eqref{4:1} or
from the TTRR \eqref{4:2}, this case is
reduced to subsection \ref{subsub411}.

\subsubsection{} \label{subsub413} $b=q^{-N}$  and
$a,c,abc^{-1}\notin \Omega(q)\setminus\{q^{-N}\}$.
 The orthogonality in this case can be obtained
taking the limit $b\to q^{-N}$.

Multiplying relation \eqref{4:4} by the factor
$(b-q^{-N})$, taking limit $b\to q^{-N}$,
and removing some non-vanishing constants one
gets, for $n\ne m$,
{\small \begin{equation}\label{4:7}
\sum_{s'=0}^{N-1}\frac{(a^{-1}cq^{s'+1},q^{s'+1};
q)_\infty~q^{-N}q^{s'}}{(cq^{s'+1};q)_\infty~(q^{-N+s'
+1};q)_{N-s'-1}}p_n(cq^{s'+1};a,q^{-N},c;q)p_m
(cq^{s'+1};a,q^{-N},c;q) =0.
\end{equation}}
The others terms of the two series in the series
representation for the Jackson $q$-integral,
vanish after taking the limit since these series
converges uniformly for $b$ in a compact
neighborhood of $q^{-N}$.

Reversing the summation and using the identity
$$
(\alpha;q)_s=(\alpha^{-1}q^{1-s};q)_s(-\alpha)^s
q^{\binom{s}{2}},
$$
orthogonality property \eqref{4:7} can be rewritten
as
{\small
\begin{equation}\label{4:8}
\sum_{s=0}^{N-1}\frac{(ac^{-1}q^{-N+1},q^{-N+1}
;q)_s}{(c^{-1}q^{-N+1},q;q)_s}\frac{q^{(N-1)s}}
{a^{s}}
p_n(cq^{N}q^{-s};a,q^{-N},c;q)p_m(cq^{N}
q^{-s};a,q^{-N},c;q)=0.
\end{equation}}
Comparing \eqref{4:5} and \eqref{4:8},
we get
\begin{align}
p_n(x;a,q^{-N},c;q)&=c^nq^{nN}h_n^{(ac^{-1}q^{-N},c)}(c^{-1}q^{-N}x;
N-1;q) \nonumber\\&=c^nq^{nN}p_n
(c^{-1}q^{-N}x;ac^{-1}q^{-N},c,q^{-N};q).
\label{4:9}
\end{align}
The used identities are not valid for several
configurations of the parameters, however
\eqref{4:9} is also valid for these
configurations by using analytic continuation.
Thus the case treated in this subsection can be
reduced to the case considered in subsection
\S \ref{subsub411} by setting $x\mapsto c^{-1}q^{-N}x$.

It is curious that identity \eqref{4:9}
has the hypergeometric form
\begin{equation*}\begin{split}
&{}_3\varphi_2\left.\left(\begin{array}{c} q^{-n},aq^{n-N+1},
x\\ aq, cq \end{array} \right|q;q\right)\\ &=\frac{c^n
q^{nN}(ac^{-1}q^{-N+1},q^{-N+1};q)_n}
{(aq,cq;q)_n}{}_3\varphi_2\left.\left(\begin{array}{c}
q^{-n}, aq^{n-N+1}, c^{-1}q^{-N}x \\ac^{-1}
q^{-N+1}, q^{-N+1} \end{array} \right| q;q\right),
\end{split} \end{equation*}
which coincides with \cite[(3.2.6)]{gara}
$$\begin{array}{c}
\displaystyle \left.{}_3\varphi_2\left(\!\!\!\begin{array}{c}
q^{-n},\widehat{a}q^{n}, \widehat{b}\\ \widehat{d},
\widehat{e} \end{array}  \!\right| q;\frac{\widehat{d}
\widehat{e}}{\widehat{a}\widehat{b}}\right)
\!=\!\frac{(\widehat{a}q/\widehat{d},\widehat{a}
q/\widehat{e};q)_n}{(\widehat{d},\widehat{e}
;q)_n} \left(\frac{\widehat{d}\widehat{e}}
{\widehat{a}q}\right)^n \!\!\!\left.{}_3\varphi_2
\left(\!\!\!\begin{array}{c} q^{-n}, \widehat{a}q^{n},
\widehat{a}\widehat{b}q/\widehat{d}
\widehat{e}\\ \widehat{a}q/\widehat{d},
\widehat{a}q/\widehat{e} \end{array}
\!\right| q;\frac{q}{\widehat{b}}\right)
\end{array}  $$
in the parameters but it does not in the
arguments if one sets $\widehat{a}=aq^{-N+1}$,
$\widehat{b}=x$, $\widehat{d}=aq$, and
$\widehat{e}=cq$.

\subsubsection{}\label{subsub414}
$abc^{-1}=q^{-N}$  and $a,b,c\notin\Omega(q)
\setminus\{q^{-N}\}$.
Once again, by \eqref{4:6}, this case can
be reduced to the case in subsection \S
\ref{subsub413}.

%%%%%%%%%%%%%%%%%%%%%%%%%%%%%%%%%%%%%%%%%%%%%%%%%
\subsection{The orthogonality conditions for $|q|\ge1$.}
\label{subsec4.2}

Identities \eqref{3:6} and (3.2.2) in \cite{gara}
$$
{}_3\varphi_2\left.\left(\begin{array}{c} q^{-n},a,b \\
d,e \end{array}\right|q; q\right)=\frac{(e/a;q)_n}{(e;q)_n}
a^n{}_3\varphi_{2}\left.\left(\begin{array}{c} q^{-n},a,d/b
\\ d,a q^{1-n}/e \end{array} \right| q; \frac{bq}{e}\right),
$$
yield
\begin{align*}
&{}_3\varphi_2\left.\left(\begin{array}{c} q^{-n}, a b q^{n+1},
x \\aq, cq\end{array}\right| q; q\right)
\\
=&\frac{(c(abq)^{-1};q^{-1})_n}{(cq^{n};q^{-1})_n}
(ab)^{n}q^{n^2+n}{}_3\varphi_2\left.\left(\begin{array}{c}
q^{n}, (ab)^{-1}q^{-n-1}, (aq)^{-1}x \\(aq)^{-1}, c(abq)^{-1}
\end{array}\right| q^{-1}; q^{-1}\right).
\end{align*}
which in terms of
big $q$-Jacobi polynomials writes as
$$
p_n(x;a,b,c;q)=\frac{1}{(a^{-1}q^{-1})^n}\,
p_n(a^{-1}q^{-1} x;a^{-1},b^{-1},c a^{-1}
b^{-1}; q^{-1}).
$$
Hence the orthogonality conditions for
big $q$-Jacobi polynomials with $|q|>1$
follow from section \ref{sec4.1}.

If $q$ is a primitive root of unity, i.e.
$q=e^{2\pi i M/N}$ with $\gcd(N,M)=1$ then
$\{k N:k\in\mathbb{N}\}\subseteq\Lambda$, and as
we did for the Askey-Wilson polynomials, the set of
big $q$-Jacobi polynomials $(p_n(x;a,b,c;q))_{n=0}^N$
under the assumptions
$$
a, b, c, ab, abc^{-1}\neq q^{k},
\qquad k=0,\dots, N-1,
$$
are uniquely determined by the
orthogonality conditions
$$
\mathscr{L}_0(p_np_m)=d_n^2 \delta_{n,m},
\quad d_n^2\ne 0,
$$
being
$$
\mathscr{L}_0(p)=\sum_{j=0}^{N-1} \omega_0
\frac{(ra^{-1}q^j,rc^{-1}q^j;q)_\infty}
{(rq^j,rbc^{-1}q^j;q)_\infty}q^sp(rq^j),
$$
with initial condition $\mathscr{L}_0(1)=1$,
and $r$ the root with minimal argument
of the equation
$$
r^N=\frac{a^N+c^N-(ab)^N-(ac)^N}{1-(ab)^N}.
$$
Moreover, since for $n\geq N$
$$
\mathscr{D}^{N}_q p_n(x;a,b,c;q)=
\frac{(q^{n-N+1};q)_N}{(1-q)^N}p_{n-N}(x;a,b,c;q),
$$

the orthogonality conditions that characterizes
big $q$-Jacobi polynomials in such case are
$$
\langle p_n,p_m\rangle=
\displaystyle \sum_{j=0}^{\infty}
\mathscr{L}_0(\mathscr{D}_q^{Nj}(p_n)
\mathscr{D}_q^{Nj}(p_m)).
$$

\begin{remark}  \label{rem4.2}
In \cite{Zhedanov} the particular case $c=1$ is
considered and, in such a case, they got
$$
\omega_s=\frac{(1-a^N)(1-abq)(b;q)_s}
{aq(b-1)(1-a^Nb^N)(a^{-1};q)_s}q^s.
$$
\end{remark}
%%%%%%%%%%%%%%%%%%%%%%%%%%%%%%%%%%%%%%%%%%%%%%%%%
\section{Extending orthogonality properties valid
up to degree $N$}
\label{sec5}
The aims of this section are for one side
to present the factorization for those
$q$-polynomials for which there exists an $N$
such that $\gamma_N=0$, and hence an orthogonality
until degree $N$ takes place, and for the other
we extend that orthogonality properties for all
non-negative degrees obtaining a Sobolev
type orthogonality properties.

Taking into account the basic idea about how
the factorization process works is already
known (see e.g. \cite{cola1}) we only show the
sketch regarding the factorization for the
$q$-polynomials.

Since $q$-polynomials fulfill, for $n\ge 0$, the TTRR 
$$
p_{n+1}(x)=(x-\beta_n)p_n(x)-\gamma_{n}p_{n-1}(x),
$$
with $p_{-1}\equiv 0$, $p_0(x)=1$, observe that if 
there exists some integer $N>0$ so that $\gamma_N=0$, 
then it is straightforward to check that, for $n\ge N$,
the following relation holds:
\begin{equation}\label{5:1}
p_n=p_N p_{n-N}^{(N)},
\end{equation}
where $(p_n^{(N)})$ is the family of $N$th associated
polynomials which fulfills, for $n\ge 0$, the recurrence
relation:
$$
p_{n+1}^{(N)}(x)=(x-\beta_{n+N})p_n^{(N)}(x)-
\gamma_{n+N}p_{n-1}^{(N)}(x),
$$
with initial conditions
$p_{-1}^{(N)}(x)\equiv0$, $p_0^{(N)}(x)=1$.

Notice that in the case of $q$-polynomials the existence
of an integer $N$ so that $\gamma_N=0$ is directly related
with the fact that there is a term of the form $q^{-N+1}$
in the denominator parameters of one of the hypergeometric
representations (see \eqref{3:1} and \eqref{4:2}).
In such a case the hypergeometric function
$_p\varphi_{p-1}$ with a suitable normalization factorizes as follows:

Let $a=\{a_1,\dots,a_{p-1}\}$ and $b=\{b_1,\dots,b_{p-2}\}$, then
\begin{align*}
&(q^{-N+1};q)_{n+N}{}_p\varphi_{p-1}\left.\left(
\begin{array}{c}q^{-n-N},a\\ q^{-N+1},b\end{array}
\right| q; z\right)\\ =&\frac{(q^{-n-N};q)_{N}\,
(a;q)_{N}z^{N}}{(b;q)_{N}(q;q)_{N}}\sum_{k=0}^{n}
\frac{(q^{-n};q)_{k}(aq^{N};q)_{k} z^k}{(bq^{N};q)_k
(q^{N+1};q)_k}~\frac{(q;q)_{n}}{(q;q)_k}\\
=&\frac{(q^{N+1};q)_{n}(a;q)_{N}z^{N}}
{(b;q)_{N}}(-1)^{N}q^{N(-n-N)+N(N-1)/2}\,
{}_p\varphi_{p-1}\left.\left(\begin{array}{c}
q^{-n},aq^{N}\\ q^{N+1}, bq^{N}\end{array} \right|
q; z\right).
\end{align*}
Hence it is straightforward combining to obtain
the following factorization:
\begin{equation}\label{5:2}
\begin{split}
(q^{-N+1};q)_{n+N}{}_p\varphi_{p-1}&\left.\left(
\begin{array}{c} q^{-n-N},a\\ q^{-N+1},b\end{array}
\right|q;z\right)=q^{-nN}(q^{N+1};q)_n(q^{-N+1};q)_N\\
&\times {}_p\varphi_{p-1}\left.\left(
\begin{array}{c} q^{-N},a\\ q^{-N+1},b\end{array}
\right|q;z\right)
{}_p\varphi_{p-1}\left.\left(
\begin{array}{c} q^{-n},aq^N\\ q^{N+1},bq^N\end{array}
\right|q;z\right).
\end{split}\end{equation}
Notice that the first hypergeometric function of the
right-hand side of \eqref{5:2}, with its
corresponding normalization coefficient, is the
polynomial of degree $N$ and the second one is
the $N$th associated polynomial in the factorization
\eqref{5:1} for $n\to n+N$.
Table \ref{table4} shows the $N$th associated
polynomial.

In the sequel we are going to assume that
no element of $b$ belongs to
$\Omega(q)$, therefore theorem \ref{theo2.6}
is applicable in such a case.
Let us go on to describe how to obtain
functionals $\mathscr L_0$, $\mathscr L_N$
and the linear operator ${\mathscr T}^{(N)}=
{\mathscr T}^{N}$ in \eqref{2:3}.

Obviously $\mathscr L_0(p)\stackrel{\rm def}=
\langle {\bf u}, p \rangle$ where $\bf u$ is
the linear form with respect to the co\-rres\-pon\-ding
family of $q$-polynomials $(p_n)_{n=1}^N$ is
orthogonal. Moreover, due the difference 
pro\-per\-ti\-es of such families $\mathscr T$ 
is going to be a difference operator and
$\mathscr L_N(p)\stackrel{\rm def}= \langle
{\bf v}, p \rangle$ where $\bf v$ is the linear
form with respect to the polynomial sequence
$({\mathscr T}^{(N)} p_{n+N})$ is orthogonal
\cite{coma2}.

Let us describe briefly the most complicated case:
the $q$-Racah polynomials.

Notice that setting $\alpha=q^{-N}$ then the $N$th
$\gamma$'s coefficient for $q$-Racah polynomial
vanishes \cite[(3.2.3)]{kost}, i.e. $\gamma_N=0$,
and therefore we can apply theorem \ref{theo2.6}, obtaining:
$$
r_n^{(N)}(x;q^{N},\beta,\gamma,\delta;q)=p_n(\mu(x)/(2
\sqrt{\gamma \delta q});q^{N}\sqrt{\gamma\delta q},
\sqrt{q/\gamma\delta},\beta \sqrt{\delta q/\gamma},
\sqrt{\gamma q/\delta};q),
$$
Moreover, taking into account that for these polynomials
$$
\frac{\Delta}{\Delta \mu(x)}r_n(\mu(x);\alpha,\beta,
\gamma,\delta|q)=\frac{q^{-n}-1}{q^{-1}-1}
r_{n-1}(\mu(x);\alpha q\beta q,\gamma q,\delta|q),
$$
and their connection with the Askey-Wilson
polynomials (see table \ref{table3}) it is
clear that the operator
$\mathscr T=\Delta/\Delta \mu(x)$ for which
we obtain that $\mathscr T^{N}\left(r_n(x;\alpha,
\beta,\gamma, \delta;q) \right)$ is, up to a constant,
equal to
$$
p_n(\mu(x)/(2\sqrt{\gamma \delta q^{N+1}});
\sqrt{\gamma\delta q^{N+1}},\sqrt{q/\gamma\delta q^N},
\beta \sqrt{\delta q^{N+1}/\gamma},\sqrt{\gamma
q^{N+1}/\delta};q).
$$
Thus the linear functional $\bf v$ is related with the
linear operator of Askey-Wilson polynomials with
parameters $\sqrt{\gamma\delta q^{N+1}}$,
$\sqrt{q/\gamma\delta q^N}$, $\beta \sqrt{\delta q^{N+1}/\gamma}$,
and $\sqrt{\gamma q^{N+1}/\delta}$.
%%%%%%%%%%%%%%%%%%%%%%%%%%%%%%%%%
\begin{table}[!hbt]
\begin{tabular}{llll}
\hline \multicolumn{2}{c}{Family and condition}&
\multicolumn{2}{c}{$N$th associated polynomial}\\ $q$H&$N\mapsto N-1$
&b$q$J&$p_n(xq^N;\alpha q^N, \beta q^N,q^N;q)$\\
d$q$H&$N\mapsto N-1$& cd$q$H&$p_n(\mu(x)/
(2\sqrt{\gamma\delta q});q^N \sqrt{\gamma\delta q},
\sqrt{q/\gamma\delta}, \sqrt{\gamma q/\delta}|q)$\\
$q$K&$N\mapsto N-1$ &b$q$J&$p_n(xq^N;q^N,-pq^{N-1},0;q)$\\
Q$q$K&$N\mapsto N-1$ &$q$M&$m_n(xq^{N};q^{N},-p^{-1}
q^{-N};q)$\\ A$q$K&$N\mapsto N-1$& b$q$L&$p_n(xq^{N};pq^{N},q^{N};q)$\\
d$q$K & $N\mapsto N-1$ &cd$q$H&$p_n(\lambda(x)/(2\sqrt{cq^{1-N}});\sqrt{cq^{N+1}},
\sqrt{q^{N+1}/c},0|q)$\\
\hline
\end{tabular}
\caption{$N$th associated polynomials involved in
the factorization \eqref{5:1}}\label{table4}
\end{table}
%%%%%%%%%%%%%%%%%%%%%%%%%%%%%%%%%%%%%%%%%%%%%%%%%%%%%%%%%%%%
%%%%%%%%%%%%%%%%%%%%%%%%%%%%%%%%%%%%%%%%%%%%%%%%%%%%%%%%%%%%
\begin{table}[!hbt]
\begin{tabular}{lll}
\hline $p_n$& ${\mathscr T}$ &${\mathscr T}^N(p_{n+N})$\\
\hline $q$H& $\frac{\Delta}{\Delta q^{-x}}$
&$p_n(x;\alpha q^N, \beta q^N,1;q)$\\ d$q$H
&$ \frac{\Delta}{\Delta \mu(x)}$ &$p_n(\mu(x)/(2\sqrt{\gamma\delta
q^{N+1}});\sqrt{\gamma\delta q^{N+1}}, \sqrt{q/\gamma\delta q^N},
\sqrt{\gamma q^{N+1}/\delta}|q)$\\
$q$K & $\frac{\Delta}{\Delta q^{-x}}$ &$p_n(x;1,-pq^{2N-1},0;q)$\\
Q$q$K &$\frac{\Delta}{\Delta q^{-x}}$ & $m_n(x;1,-p^{-1}q^{-N};q)$\\
A$q$K &$\frac{\Delta}{\Delta q^{-x}}$ & $p_n(x;pq^{N},1;q)$\\
d$q$K & $\frac{\Delta}{\Delta \lambda(x)}$&
$p_n(\lambda(x)/(2\sqrt{cq});\sqrt{cq}, \sqrt{q/c},0|q)$\\
\hline
\end{tabular}
\caption{Unnormalized ${\mathscr T}^{(N)}(p_{n+N})$
involved in factorization \eqref{5:1}}\label{table5}
\end{table}
%%%%%%%%%%%%%%%%%%%%%%%%%%%%%%%%%%%%%%%%%%%%%%%%%%%%%%%%%%%%

%%%%%%%%%%%%%%%%%%%%%%%%%%%%%%%%%%%%%%%%%%%%%%%%%%%%%%%%%%%%
\begin{table}[!hbt]
\begin{tabular}{lll}
\hline
$q$R$\to$ AW&
$r_n(x;\alpha,\beta,\gamma,\delta;q)$ &
$p_n\left(\frac{x}{2\sqrt{\gamma\delta q}};
\sqrt{\gamma\delta q},\,\alpha\sqrt{
\frac{q}{\gamma\delta}},\, \beta\sqrt{
\frac{\delta q}{\gamma}},\,\sqrt{\frac{\gamma
q}{\delta}};q\right)$\\
AW$\to q$R&
$p_n(x;a,b,c,d;q)$& $r_n\left(2ax;\frac{ab}{q},
\frac{cd}{q},\frac{ad}{q},\frac{a}{d};q\right)$\\
b$q$J$\to q$H&$p_n(x;a,b,c;q)$ &$h_n(x;a,b,-1-\log_q
c\,;q)$\\ $q$H$\to$b$q$J&$h_n(x;a,b,N;q)$&$p_n(x;a,b,
q^{-N-1};q)$\\d$q$H$\to$cd$q$H&$r_n(x;\gamma,\delta,
N;q)$&$p_n\left(\frac{x}{2\sqrt{\gamma\delta q}};
\sqrt{\gamma\delta q},\,\sqrt{\frac{\gamma q}
{\delta}},\,\frac{1}{q^N\sqrt{\gamma\delta q}}
|q\right)$\\cd$q$H$\to$d$q$H&$p_n(x;a,b,c|q)$&$ r_n
\left(2ax;\frac{ab}{q},\frac{a}{b},-\log_q(ac);
q\right)$\\Q$q$K$\to q$M&$k^{qtm}_n(x;p,N;q)$&$
m_n\left(x;q^{-N-1},-\frac{1}{p};q\right)$\\
$q$M$\to$Q$q$K&$m_n(x;b,c;q)$& $ k^{qtm}_n\left(x;
-1-\log_q b,-\frac{1}{c};q\right)$\\ Q$q$K$\to$A$q$K
&$k^{qtm}_n(x;p,N;q)$&$ k^{aff}_n\left(x q^{N};
p^{-1},N;q^{-1}\right)$\\A$q$K$\to$Q$q$K&$k^{aff}_n
(x;p,N;q)$&$k^{qtm}_n(xq^{-N};p^{-1},N;q^{-1})$
\\ $q$K$\to$l$q$J&$k_n(x;p,N;q)$&$
p_n(x q^N;-pq^N,q^{-N-1};q)$\\ l$q$J$\to q$K
&$p_n(x;a,b;q)$&$k_n(bqx;-abq,-1-\log_qb;q)$\\
A$q$K$\to$b$q$L&$k^{aff}_n(x;p,N;q)$&$p_n(x;p,q^{-N-1};
q)$\\b$q$L$\to$A$q$K&$p_n(x;a,b;q)$&$k^{aff}_n
(x;a,-1-\log_qN;q)$ \\
l$q$J$\to$b$q$J&$p_n(x;a,b;q)$& $p_n(bqx; b,a,0;q)$\\
$q$K$\to$ b$q$J&$k_n(x;p,N;q)$& $p_n(x;q^{-N-1},-pq^N,0;q)$
\\ \hline
\end{tabular}
\caption{Some unnormalized identities between
$q$-polynomials.}\label{table3}
\end{table}
%%%%%%%%%%%%%%%%%%%%%%%%%%%%%%%%%%%%%%%%%%%%%%%%%%%%%%%%%%%%
%%%%%%%%%%%%%%%%%%%%%%%%%%%%%%%%%%%%%%%%%%%%%%%%%%%%%%%%%%%%
\bibliographystyle{plain}

\end{document}